
\magnification=\magstep1

\input amstex

\documentstyle{amsppt}
\input epsf.tex

\TagsOnRight

\def\cite#1{{\rm [#1]}}

\leftheadtext{P. ERD\H{O}S, E. MAKAI, JR., J. PACH}

\rightheadtext{Two nearly equal distances in ${\Bbb{R}}^d$}

\topmatter

\title
Two nearly equal distances in $R^d$
\endtitle

\author
Paul Erd\H{o}s, Endre Makai, Jr., J\'anos Pach
\vskip.5cm
\eightpoint
Alfr\'ed R\'enyi Institute of Mathematics, HUN-REN,
1364 Budapest, POB 127, Hungary$^{*}$\\
\vskip.1cm
{\lowercase{\rm{makai.endre}}}{\@}{\lowercase{\rm{renyi.hu,}}}
{\lowercase{\rm{http://www.renyi.hu/\~{}makai}}}\\
{\rm{ORCID ID: https://orcid.org/0000-0002-1423-8613}}\\
{\lowercase{\rm{pach}}}{\@}{\lowercase{\rm{renyi.hu,}}}
{\lowercase{\rm{http://www.renyi.hu/\~{}pach}}}\\
{\rm{ORCID ID: https://orcid.org/0000-0002-2389-2035}}
\vskip.1cm
{\it 2020 Mathematics Subject Classification:} Primary:52C10,
Secondary: 05C35
\vskip.1cm
{\it Keywords:}
{\rm  Erd\H{o}s-type problems, separated point sets,
nearly equal distances, $k$-distance sets, Tur\'an's theorem,
Ramsey's theorem}
\endauthor


\abstract
A point set $P \subset {\Bbb{R}}^d$ is {\it separated} if the
minimum distance between any two points in $P$ is at least~$1$.
For $d \ne 4,5,$ we determine, for every $t_1,t_2 \ge 1$,
and for $n$ at least a suitable $n_d$, the maximum number of
point pairs in a separated $n$-element point set in ${\Bbb{R}}^d$,
with distances in the set $[t_1,t_1 + 1]\cup[t_2,t_2 + 1]$.
For $d=4,5$ we establish a weaker, similar asymptotic estimate.
Recently N. Frankl and A. Kupavskii have generalized this result
to unions of $k\ge 2$ intervals. We also determine the maximum
number of point pairs in an $n$-element point set in ${\Bbb{R}}^d$,
whose distances belong to the union of $k \ge 2$ intervals of the
form $[t_i, t_i(1 + \varepsilon)]$, where $t_i > 0$ and 
$\varepsilon > 0$ is small.
\endabstract

\endtopmatter

\document


\heading
\S1. Introduction
\endheading

Around 1945, Paul Erd\H os found two interesting applications of extremal combinatorics. One is related to an algebraic question of Littlewood and Offord~\cite
{19}, 
 and the other one is
in geometry. In~\cite
{6}, 
 he applied Sperner's lemma to give a tight upper bound on the number of subsets of a set of $n$ real numbers, whose absolute values are at least $1$, and whose sums fall into a given interval of length $1$. In~\cite
{7}, 
 Erd\H os addressed the following question: At most how many times can the same distance occur among $n$ points in ${\Bbb{R}}^d$? More precisely, what is the maximum number of unordered point pairs
that determine the same distance?

Erd\H os
modified the second question in the spirit of the first one,
cf. \cite
{12}. 
 At most how many unordered pairs $\{p,q\}$ of distinct points can be selected from an $n$-element point set $P \subset
{\Bbb{R}}^d$ so that all distances $d(p,q)$ are {\it nearly the same}, in the sense that they fall into the same unit interval? To avoid the degenerate situation where all points are very close to
each other and, hence, all distances are nearly $0$, we consider only {\it separated} point sets $P$. That is, we assume that the distance between any two points of $P$ is at least $1$.
To give an answer to the last question, we recall {Tur\'an's
theorem}~\cite
{26, 1}. 
 For $n,k \ge 1$ integers, define the {\it{Tur\'an number}} $T(k,n)$,
as the maximum number of edges that a graph on $n$ vertices can have without
containing a complete subgraph $K_k$ on $k$ vertices. According to {\it{Tur\'an's theorem}}, for a fixed $k$, we have
$$
T(k,n) = {n^2\over 2}\left( 1 - {1\over k-1} \right) + O_k(1) \le
{n^2\over 2}\left( 1 - {1\over k-1} \right).
$$
Moreover, the only $K_k$-free graph for which this maximum is attained is the so-called Tur\'an graph. This is a complete
$(k-1)$-partite graph whose classes are as equal as possible, i.e., each class consists of $\lfloor n / (k-1) \rfloor $ or
$\lceil n / (k-1) \rceil $ points.
\medskip

{\bf{Theorem A.}} (\cite
{12}, 
 Theorem 5)
{\it{For any $d\ge 2$, there exist positive constants $c_d, n_d$ such that for every $t\ge 1$, every separated set $P\subset {\Bbb{R}}^d$ with $n\ge n_d$ elements the following holds. The
set $P$
has at most $T(d+1,n)$ unordered
point pairs whose distances belong to the interval
$[t,t+c_d n^{1/d}]$. This bound is best possible for every
$d$ and every $n \ge n_d$.}}
\medskip

To see that the bound $T(d+1,n)$ can be attained, we write ${\Bbb{R}}^d$ as
${\Bbb{R}}^{d-1} \times {\Bbb{R}}$, and let $q_1, \ldots ,
q_d \in {\Bbb{R}}^{d-1}$ be the vertices of a regular $(d-1)$-simplex
of edge length $t$. Write $n$ as a sum,
$n=n_1+\ldots+n_d$, where $n_i = \lfloor n/d \rfloor $ or
$\lceil n/d \rceil $ for every $i$.
Let $P := \{ q_i + j e_d : 1 \le i \le d, \,\, 1 \le j \le n_i
\} $, where $e_d = (0, \ldots , 0,1)$.
If $t$ is large enough (depending on $n$), then all distances between two points in distinct sets $P_{i(1)}$ and $P_{i(2)}$
belong to the interval $[t,t+1]$, and the number of such pairs is $T(d+1,n)$.

Originally, Theorem A was stated for unit intervals $[t,t+1]$, but
its proof easily extends to this case. (See the paragraph after Lemma 3.1 in \cite
{12}.) 

\medskip

We say that a set {\it determines a distance} $t>0$ if it
has two points at distance $t$ from each other.
It is our goal to extend Theorem A and obtain an upper bound for the number of pairs
whose distances fall into the union of $k\ge 2$ unit, or short, intervals.
In \cite
{9}, 
 we
made the first step in this direction by providing an asymptotically tight bound in the plane.

\medskip

{\bf{Theorem B.}} (\cite
{9}, 
 Theorem 2)
{\it{For any $k\ge 2$ and $\varepsilon>0$, there exist positive constants $c_ {k,\varepsilon}$ and $n_{k, \varepsilon }$ such that for every $t_1,\ldots,t_k\ge 1$, for every separated set $P\subset{\Bbb{R}}^2$ with $|P|=n\ge n_{k,\varepsilon}$,
the following holds.

The number of unordered
point pairs from $P$ that determine a distance belonging to the set
\,
$\cup _{i=1}^k [t_i,t_i + c_{k, \varepsilon } n ^{1/2}]$,
is at most
$$
\frac{n^2}{2}\left(1 - \frac{1}{k+1} + \varepsilon \right).
$$
This statement is asymptotically tight: it does not remain true if we replace  the last expression by $T(k+2,n) - 1$.}}

\medskip

An example of an $n$-element point set with $T(k+2,n)$ pairs whose distances are nearly equal to one of $k$ numbers,
$t_1,\ldots,t_k$ is the following. Let $t_i :=it$, for $1 \le i \le k$, and let $n=n_1+\ldots+n_{k+1},$ where the $n_h$'s, for $1 \le h \le k+1$, are as equal as possible.
Let
$P_h=\{((h-1)t, j) : 1\le j\le n_h\}$ and $P=\cup_{h=1}^{k+1}P_h$.
If, for a given $n$, $t$ is large enough, then every
distance between two points belonging to distinct $P_h$'s lies in $\cup_{i=1}^k[t_i,t_i+1]$.
\medskip

To generalize Theorem B to higher dimensions, we need a definition.

\medskip

{\bf{Definition 1.}}
For any positive integers $d$ and $k$, we call a finite subset of ${\Bbb{R}}^d$ a {\it $k$-distance set} if it determines at most $k$ distinct (positive) distances. Let $m(d,k)$
denote the maximum cardinality of a $k$-distance set in ${\Bbb{R}}^d$. (This exists by Ramsey's theorem.)
If $k=2$,
we write $m(d) := m(d,2)$, for simplicity.
\medskip

Estimating the value of $m(d,k)$ is equivalent to Erd\H os's {\it distinct distances problem} \cite
{6, 7} 
 and has a huge literature. In particular, it is known \cite
{2, 3} 
 that
$$
\binom{d+1}{k} \le m(d,k) \le \binom{d+k}{k}.
\tag 1.1
$$
This implies that
for a fixed $k$ and $d \to \infty $, we have
$m(d,k) = (d^k / k!) \left( 1 + o_k(1) \right) $, while
for a fixed $d$ and $k \to \infty $, we have
$m(d,k) \le (k^d / d!) \left( 1 + o_d(1) \right) $.
The asymptotically best upper bounds for $m(2,k)$ and $m(d,k)$ for
$d \ge 3$
have been established by Guth and Katz \cite
{16} 
 and by Solymosi and Vu \cite
{25}, 
resp.

For our purposes, the case $k=2$ will be relevant. For the maximum cardinality $m(d)=m(d,2)$ of a $2$-distance set in ${\Bbb{R}}^d$, it is known that
$$
\cases
m(1) = 3,\,\, m(2) = 5\,\,  {\text{\cite
{8}}}, 
 \,\, m(3) = 6\,\,
{\text{\cite
{5}}}, 
 \,\,m(4) = 10,
\\
m(5) = 16,\,\, m(6) = 27,\,\, m(7) = 29,\,\,
m(8) = 45\,\, {\text{\cite
{18}}}. 
\endcases
\tag 1.2
$$
In particular,
$$
{\text{for all }} d \ge 2, {\text{ we have }} m(d-1) > d \,.
\tag 1.3
$$
\medskip

Our main result is the following generalization of the
special case $k=2$ of Theorem B to higher dimensions.
\medskip

{\bf{Theorem 1.}}
{\it{For any integer $d \geq 2$, $d\not=4,5$, there exist positive constants $c_d, n_d$ such that for any $t_1,t_2 \ge 1$, for
every separated point set $P\subset{\Bbb{R}}^d$ with
$n \ge n_d$ elements, the following holds. 
The number of unordered point pairs in $P$ that determine a distance belonging to the set $[t_1,t_1 + c_d
n^{1/d}]\cup [t_2,t_2 + c_d n^{1/d}]$,
is at most
$$
T\left( m(d - 1) + 1, n \right) = {n^2\over 2}\left(1 - {1\over m(d - 1)}\right)+ O_d(1).
$$

For $d=4$ or $5$, for any $\varepsilon > 0$, there exist positive constants $c_{d, \varepsilon }, n_{d, \varepsilon }$ such that for any $t_1,t_2 \ge 1$, for every separated point set $P\subset{\Bbb{R}}^d$ with $n\ge n_{d, \varepsilon }$ elements,
the following holds. The
number of unordered point pairs in $P$
that determine a distance belonging to the set $[t_1,t_1 + c_{d, \varepsilon }(\log n)^{1/d}]\cup
[t_2,t_2 + c_{d,\varepsilon }(\log n)^{1/d}]$ is at most
$$
{n^2\over 2}\left(1 - {1\over m(d - 1)} + \varepsilon \right).
$$

These upper bounds cannot be reduced to $T\left( m(d - 1) + 1,
n \right) - 1$,
for any $d$ and $n$.}}

\medskip

We also study a closely related problem, where two distances are considered nearly equal if they fall into an interval $[t,t(1+\varepsilon)]$, for some small $\varepsilon>0$. To formulate our result we need to extend Definition 1, as follows.
\medskip

{\bf{Definition 2.}}
For any $\varepsilon \ge 0$ and integers $d,k \ge 1$, we call a finite
subset of ${\Bbb{R}}^d$ a {\it{$(k, \varepsilon )$-distance set}} if all distances
determined by it lie in the union of $k$ intervals of the form
$[t_1, t_1(1+ \varepsilon)], \ldots, [t_k, t_k(1 + \varepsilon)]$,
for some $t_1, \ldots,t_k > 0$. Let $m(d,k,\varepsilon)$ denote the maximal
cardinality of a $(k, \varepsilon )$-distance set in ${\Bbb{R}}^d$. (This is finite for every $\varepsilon > 0$.
In fact, by
applying Ramsey's theorem, it is
enough to see that $m(d,1,\varepsilon)$ is finite, and this
follows from a volume argument.) 

Obviously, a $(k,0)$-distance set is a $k$-distance set and
$m(d,k,0) = m(d,k)$.


\medskip

{\bf{Theorem 2.}}
{\it{For any fixed integers $d, k \geq 1$ there exists $\varepsilon_{d,k} > 0$
such that for $0 < \varepsilon < \varepsilon_{d,k}$ the following two
statements hold.

{\rm (A)} For the maximum cardinality of a $(k, \varepsilon )$-distance set
in ${\Bbb{R}}^d$, we have
$$
m(d,k,\varepsilon)=(d + 1)^k.
$$

{\rm (B)} For any set $P\subset{\Bbb{R}}^d$
of \,$n\ge 1$ points, and for any \,$t_1,\ldots,t_k>0$, the following holds.
The number of unordered pairs in $P$
that determine a distance belonging to the set\,\,\,\, $\cup_{j = 1}^k$ $[t_j, t_j(1 + \varepsilon)]$,
is at most
the Tur\'an number $T \left( (d + 1)^k + 1, n \right) $.
This upper bound cannot be reduced to
$T\left( (d + 1)^k + 1 , n \right) - 1$, for any $d,k$ and $n$, and any
$\varepsilon > 0$.}}

\medskip


It follows from Theorem 2 (A) and \thetag{1.1}
that, for $k$ fixed and $d \to \infty $,
$$
1\ge \frac{m(d,k)}{\lim\limits_{\varepsilon\searrow 0} m(d,k,
\varepsilon)}
=\frac{m(d,k)}{(d+1)^k} = \frac{1}{k!}+o_k(1).
$$


Observe that in Definition 2 and Theorem 2, the assumption that
$P$ is separated is not required. (Actually,
the concept of a $(k, \varepsilon )$-distance
set is similarity invariant, so 
we could have required this property as well.)


\medskip

The rest of this paper is organized as follows. In \S2, we describe several constructions showing the tightness of Theorems 1 and 2. \S3 and \S4 contain the proofs of Theorems 1 and 2, resp. In
\S5, we make some concluding remarks.
\medskip


The present paper is a minimally edited version of a manuscript written in the early 1990s.
We posted it on arXiv in January 2019 \cite
{10}. 
A somewhat weaker version of Theorem 1 was announced in \cite
{21} 
 in 2002.
Our proofs use simple Tur\'an-type results and elementary geometric observations. The first inequality of
Theorem 1 has been generalized by N\'ora Frankl and Andrey
Kupavskii 
to unions of $k \ge 2$ intervals, for
any $d \ge d(k)$ for some $d(k)$
(\cite
{14} 
 Theorem 1.2 and
\cite
{15} 
 Theorem 13).

Moreover, they proved in \cite
{15} 
 Theorem 12 the following. Let us
fix any
$d,k \ge 2$. Then there exists a natural number $N_k(d)$, such
that the following holds. For any $\varepsilon > 0$, there exists
a natural number $n(d,k,\varepsilon )$, such that for all $n \ge
n(d,k,\varepsilon )$ the following is valid. 
The maximum 
number of unordered pairs of points, from any $n$
points in ${\Bbb{R}}^d$, whose distances lie in the
union of $k$ intervals,
lies in $[T(N_k(d),n), T(N_k(d),n) + \varepsilon n^2]$.

\medskip


\heading
\S2. Constructions
\endheading

The aim of this section is to describe the constructions showing the tightness of
Theorems 1 and 2.
\medskip

{\bf{Construction 1.}}
We regard ${\Bbb{R}}^{d - 1}$ as the hyperplane of ${\Bbb{R}}^d$
spanned by the first $d-1$ usual basic unit vectors.
Let $Q \subset {\Bbb{R}}^{d - 1}$
be a finite point set, with all distances sufficiently large, and let $m :=|Q|$.
Suppose that the
distances determined by $Q$ all lie in the union of $k$
intervals of length $\varepsilon $ each, where $0 \leq \varepsilon < 1$.
Let $Q = \{q_1, \ldots, q_m\}$.
Let $n = n_1 + \ldots + n_m$, where each $n_i$ is $\lfloor n/m
\rfloor $ or $\lceil n/m \rceil $.
We construct a point system $P = P(Q)$
of $n$ points in ${\Bbb{R}}^d = {\Bbb{R}}^{d - 1} \times {\Bbb{R}}$ as
follows.
We let
$$
P(Q) := \{ q_i + j e_d: 1 \leq i \leq m, \,\, 1 \leq j \leq n_i \} ,
$$
where $e_d$ is the
$d$-th usual unit basic vector in ${\Bbb{R}}^d$.
If, for given $n$,
all distances determined by $Q$ are large enough, then
the following holds.
The
distances of all pairs of points $ q_{i(1)} + j(1) e_d, \,\,
q_{i(2)} + j(2) e_d \in P(Q)$
with $i(1) \neq i(2)$ lie in the union of $k$ unit intervals
(or of $k$
arbitrarily small intervals, provided $\varepsilon $ can be made arbitrarily
small).
The number of these pairs of points is $(n^2/2) (1 - 1/m) +
O_{d,k}(1) \le (n^2/2) (1 - 1/m)$, for $n \to \infty $.

\medskip


We present two particular cases of Construction 1.


\medskip

{\bf{Construction 1$'$.}} The case $k=2$ of this construction will show the
tightness of Theorem 1.

Let $k$ be fixed and $d \to \infty $.
In Construction 1, we choose
$Q \subset {\Bbb{R}}^{d-1}$ as a $k$-distance subset of
maximum cardinality $m(d-1,k)$, with all distances sufficiently large.
By \thetag{1.1},
$$
|Q| = m(d-1,k) = \frac{d^k}{k!} \left( 1 + o_k (1) \right) .
$$
Then the set $P(Q)$ determines
$$
\frac{n^2}{2}\left(1 - \frac{1}{|Q|}\right) + O_{d,k}(1)
=
\frac{n^2}{2}\left(1 - \frac{1}{m(d-1,k)}\right) + O_{d,k}(1)
$$
\vskip-0.4cm
$$
\le
\frac{n^2}{2}\left(1 - \frac{1}{m(d-1,k)}\right)
$$
distances, taken with multiplicity,
that lie in the union of $k$ intervals of arbitrarily small length.

\medskip


{\bf{Construction 1$''$.}}
Let $d$ be fixed and $k \to \infty $.
We construct a set $Q \subset {\Bbb{R}}^{d-1}$ as follows.
Let $k = k_1 + \ldots + k_{d-1}$, where each $k_i$ is $\lfloor
k/(d-1) \rfloor $ or $\lceil k/(d-1) \rceil $. We write
$\{ e_1, \ldots , e_d \} $ for the usual basic unit vectors in ${\Bbb{R}}^d$.
Let $n \ll \lambda_1 \ll \lambda_2 \ll \ldots \ll \lambda_{d-1}$ and let
$$
Q := \bigl\{ \sum _{i=1}^{d-1} j_i \lambda _i e_i
: j_i \in \{ 0,1, \ldots , k_i \} \bigr\} .
$$
Then the distance between any two distinct
points, $\sum _{i=1}^{d-1} j_{i(1)} \lambda _i e_i , \,\,
\sum _{i=1}^{d-1} j_{i(2)} \lambda _i e_i $
\newline
$\in Q $,
is very close to one of the distances
$\lambda _i, 2 \lambda _i, \ldots , k_i \lambda _i$, where $i$ is the largest
index $\ell \in \{ 1, \ldots , d-1 \} $ such
that $j_{\ell (1)} \ne j_{\ell (2)}$.
The total number of these distances is
$k_1 + \ldots + k_{d-1} = k$, and we have
$|Q| = \prod _{i=1}^{d-1}(k_i+1)$. Hence,
for a fixed $d$ and $k \to \infty $, we have
$$
|Q| = \frac{k^{d-1}}{(d-1)^{d-1}} \left( 1 + o_d (1) \right)
\le \frac{(k + d - 1)^{d-1}}{(d-1)^{d-1}}.
$$
Using that $n \ll \lambda _1$, the number of distances determined by $P(Q)$ that lie in the union of $k$ intervals of arbitrarily small length is
$$
\frac{n^2}{2}
\left(1-\frac{1}{|Q|}\right) + O_{d,k}(1)
\le \frac{n^2}{2}
\left( 1 - \frac{(d-1)^{d-1}}{(k + d - 1)^{d-1}}\right).
$$


\medskip

It is somewhat surprising
that for a fixed $d \ge 3$ and any $\varepsilon > 0$,
a point set in ${\Bbb{R}}^{d-1}$ in which all distances are at least $1$ and belong
to $k$ intervals of length $\varepsilon $, can be much larger
than the conjectured maximum size of a point set in ${\Bbb{R}}^{d-1}$ in which every point pair determines one of $k$ specific distances. (The conjectures are $m(2,k) = \Theta
\left( k (\log k) ^{1/2} \right) $, and $d-1 \ge 3
\Longrightarrow
m(d-1,k) = \Theta _d (k^{(d-1)/2})$, see \cite
{7}.) 
 This is in sharp contrast with Theorem 1.1 in \cite
{14}, 
 and Theorem 1 in \cite
{15}, 
 stating that if $k$ is fixed, $d \ge d_k$, and $\varepsilon \in (0, \varepsilon _{d,k})$, then
these two quantities coincide.


\medskip

{\bf{Construction 2.}}
We construct, for any $d,k \ge 1$ and any $\varepsilon \in (0,1)$, a
$(k, \varepsilon )$-distance set in ${\Bbb{R}}^d$, of cardinality $(d + 1)^k$.
This will show the tightness of Theorem 2, (A).

Let us choose, for some small $\varepsilon_1 \in (0,1)$, positive
numbers $s_1, \ldots ,
s_k$, satisfying $s_i/s_{i+1} \le \varepsilon_1$ for every $i \in \{ 1,
\ldots , k-1 \} $. Fix $k$ regular simplices
centred at $0$, with circumradii $s_1, \ldots , s_k$,
and with vertices
$$
\{v_{1,i} : 1 \leq i \leq d + 1\}, \ldots, \{v_{k,i} :
1 \leq i \leq d + 1\}.
$$
Define the set of $(d + 1)^k$ vectors, 
$$
S := \{v_{1, i(1)}
+ \ldots + v_{k, i(k)} : 1 \leq i(1) \leq d + 1, \ldots,
1 \leq i(k) \leq d + 1\}.
$$
For any different $v_{1, i(1)} + \ldots + v_{k, i(k)}
, \,\, v_{1, j(1)}
+ \ldots + v_{k, j(k)} \in S$, 
let $h$ be the largest index $\ell \in \{1,
\ldots , k\}$ such that $i(\ell ) \neq j(\ell )$.
Then their distance equals 
$$
d(v_{1, i(1)} + \ldots
+ v_{h, i(h)}, v_{1, j(1)} + \ldots + v_{h, j(h)}) \in
$$
\vskip-0.4cm
$$
[d(v_{h, i(h)}, v_{h, j(h)}) - 2 s_{h - 1} -  \ldots - 2 s_1,
d(v_{h, i(h)}, v_{h, j(h)}) + 2 s_{h - 1} +  \ldots + 2 s_1] =
$$
$$
[ \left( 2(1+1/d) \right) ^{1/2} s_h - 2 s_{h - 1} -  \ldots - 2 s_1,
\left( 2(1+1/d) \right) ^{1 / 2} s_h + 2 s_{h - 1} +  \ldots + 2 s_1].
$$
If $\varepsilon_1$ is sufficiently small, then for any
$h \in \{1, \ldots, k\}$ 
the quotient of the maximum and the minimum of the last
interval lies in $[1, 1 + \varepsilon ]$.
Hence, $S$ is a $(k, \varepsilon )$-distance set, with
$$
t_h := \left( 2(1+1/d) \right) ^{1 / 2}
s_h - 2 s_{h - 1} - \ldots - 2 s_1
{\text{ for any }} h \in \{ 1, \ldots , k \} .
$$

\medskip


{\bf{Construction 3.}}
We construct, for any $d, k \ge 1$,
any $\varepsilon > 0$ and any $n$, a set
$\{ p_1, \ldots, p_n \} $ of $n$ points in ${\Bbb{R}}^d$
with the following property.
The number of point pairs determining a distance   
that belongs to $\cup_{j = 1}^k [t_j, t_j (1
+ \varepsilon)]$, for some $t_1, \ldots, t_k > 0$, is equal to $T\left( (d + 1)^k + 1 , n \right). $
This will show the tightness of Theorem 2, (B).

The points $p_1, \ldots, p_n$ are divided into
$(d + 1)^k$ classes, with
$\left\lfloor n/(d + 1)^k \right\rfloor $ or
$\left\lceil n/(d + 1)^k \right\rceil $ points in each
class, so that the distance
between any two points in different classes belongs to $\cup_{j = 1}^k
[t_j, t_j(1 + \varepsilon)]$. Each of the $(d + 1)^k$ classes of points is 
chosen in the $1$-neighbourhood of one of the $(d + 1)^k$ points of the
set $S$ as in Construction 2, where we also assume that
$1/s_1 \le \varepsilon _1$.
Like in Construction 2,
the distance between any two points in different classes belongs to the interval
$$
[ \left( 2(1+1/d) \right) ^{1/2} s_h - 2 s_{h-1} - \ldots - 2 s_1 - 2,
\left( 2(1+1/d) \right) ^{1/2} s_h + 2 s_{h-1} + \ldots + 2
s_1 + 2].
$$
Here, $h$ is the largest index $\ell \in \{1, \ldots, k\}$ such that $i(\ell ) \neq j(\ell )$,
with 
$v_{1,i(1)} + \ldots + v_{k, i(k)}$ and $v_{1, j(1)} +
\ldots + v_{k, j(k)}$ being
the elements of $S$ in Construction 2,
associated with the classes of the two points.
If $\varepsilon_1$ is sufficiently small, then
for any $h \in \{1, \ldots, k\}$, 
the quotient of the maximum and the minimum elements of the last
(displayed) interval
lies in $[1, 1 + \varepsilon ]$.
Thus we can choose
$$
t_h := \left( 2(1+1/d) \right) ^{1/2} s_h - 2 s_{h-1} -
\ldots - 2 s_1 - 2,
{\text{ for any }} h \in \{ 1, \ldots , k \} .
$$

\medskip


\heading
\S 3. Proof of Theorem 1
\endheading

We first agree on some notation and terminology. We denote the vertex set of a graph $G$ by $V(G)$. Throughout this paper, the term {\it subgraph} will always stand for {\it induced} or {\it spanned subgraph}. Let $d(p,q)$ denote the {\it distance between two points} $p,q \in {\Bbb{R}}^d$. The
{\it{norm of $p \in {\Bbb{R}}^d$}} is denoted by $\| p \| $. We write $S^{d-1}$ for the unit sphere in ${\Bbb{R}^d}$.
For any set $P \subset {\Bbb{R}}^d$, we write ${\text{diam}}\,P$,
${\text{aff}}\, P$ and ${\text{lin}}\, P$
for the {\it diameter, affine hull,} and
{\it linear hull of} $P$, resp. The {\it volume} (Lebesgue
measure) {\it{of a set in ${\Bbb{R}^d}$}}
is denoted by $V( \cdot )$, while the {\it{$(d-1)$-volume}}
is denoted by
$V_{d-1}( \cdot )$.
For $x_1, \ldots , x_d \in {\Bbb{R}}^d$, denote by
${\text{det}}\,(x_1, \ldots , x_d)$ the {\it determinant
with columns}
$x_1, \ldots , x_d$.
For any $1 \le \ell \le d + 1$ and any affinely independent vectors $x_1, \ldots , x_{\ell } \in
{\Bbb{R}}^d$, let $S(x_1, \ldots , x_{\ell })$ stand for the
{\it{$(\ell - 1)$-dimensional simplex spanned by these
vertices.}}

Throughout, we suppose that $t_1 \le \ldots \le t_k$.
The interval $[t_{\kappa }, t_{\kappa } + c_d n^{1/d}]$ will
be referred to 
as the {\it{$\kappa $-th interval.}} The symbols
${\text{const}}_d, C_d, D_d, c_d$ will
denote positive constants depending on $d$ (or on other parameters in the subscript). At different places,
${\text{const}}_d$ may stand for different constants.
We always assume that $n$ is sufficiently large in terms of
all fixed parameters.

\medskip

In the rest of this section, we present the proof of Theorem 1.
The proof falls into eight simple steps marked as
{\bf{Step 1}},
{\bf{Step 2}},
etc. For $d\not= 4,5$, we give the proof in full detail. The
treatment of the cases d=4,5 requires only minor
modifications which are described in {\bf{Step 7}} below.

\medskip

{\it{Proof of Theorem 1.}}
{\bf{Step 1}}.
The tightness of Theorem 1 was shown by Construction 1$'$.
Therefore, we only have to prove the upper bounds. 
Let $P = \{ p_1, \ldots , p_n\} $. 

\medskip


{\bf{Lemma 1.}}
{\it{It is sufficient to prove Theorem 1 under the following assumptions.


{\rm{(1)}} The intervals
$[t_1, t_1 + c_d n^{1/d}]$ and $[t_2, t_2 + c_d
n^{1/d}]$
are disjoint, and both contain at least one distance between
two points of $P$. 

{\rm{(2)}} We have $t_2 > t_1 \ge C_d n^{1/d}$, where $C_d > 1$
can be chosen arbitrarily large.

{\rm{(3)}} The ratio of any two distances that belong to the
$\kappa $-th interval ($\kappa =1,2$)
lies in $
[(1 + c_d)^{-1}, 1 + c_d]$. Hence, it lies
in
an arbitrarily small neighbourhood of 1, provided that we
choose $c_d > 0$ sufficiently small.}}

\medskip


{\it{Proof.}}
(1)
If
$[t_1, t_1 + c_d n^{1/d}]\cap[t_2, t_2 + c_dn^{1/d}]\neq\emptyset,$ 
then the length of the union of the two intervals is at most
$2c_d n^{1/d}$. Hence, if $c_d > 0$ is
sufficiently small, Theorem A yields the following.
The number of pairs $\{ p_{i(1)}, p_{i(2)} \} $
whose distances belong to the union
of the two intervals is at most
$T(d + 1, n)$.
By \thetag{1.3}, we have
$ T(d + 1, n) \le T\left( m(d - 1) + 1, n \right) $,
and Theorem 1 follows.

The same argument applies if one of the intervals does not
contain any distance $d(p_{i(1)},p_{i(2)})$.

(2) Suppose that
$
t_1 \le C_d n^{1/d}
$
for an arbitrarily large constant $C_d$.
By our assumptions, the open
balls of radius $1/2$ centred at the points  $p_i$
are disjoint.
Thus, by volume considerations,
for any fixed $p_{i(1)}$, the number of $p_{i(2)}$'s with
$d(p_{i(1)}, p_{i(2)}) \in [t_1,
t_1 + c_d n^{1/d}]$
is at most
${\text{\rm{const}}}_d \cdot
\bigl[(t_1 + c_d n^{1/d} + 1/2)^d - (t_1 - 1/2)^d \bigr]$.
Hence, the number of all pairs $\{ p_{i(1)}, p_{i(2) } \} $,
where $d(p_{i(1)},p_{i(2)})$ belongs to
the first interval, is at most
$$
{\text{\rm{const}}}_d \cdot
\bigl[n \left( (t_1 + c_d n^{1/d} + 1/2)^d -
(t_1 - 1/2)^d \right) \bigr]
\le
n \cdot {\text{\rm{const}}}_d \cdot (t_1 + c_d n^{1/d})^{d - 1}
\cdot c_d n^{1/d}
$$
\vskip-1cm
$$
\le n \cdot {\text{\rm{const}}}_d \cdot (C_d n^{1/d} + c_d n^{1/d})^{d-1}
\cdot c_d n^{1/d} 
=
n^2 \cdot {\text{\rm{const}}}_d \cdot (C_d + c_d)^{d-1}c_d \le \delta n^2,
$$
provided that we choose $c_d > 0$ so small compared to
${\text{\rm{const}}}_d$ and $C_d$ that
${\text{\rm{const}}}_d \cdot (C_d + c_d)^{d-1}c_d \le \delta $ holds.

By Theorem A, the number of pairs $\{ p_{i(1)}, p_{i(2)} \} $
with
$d(p_{i(1)},p_{i(2)}) \in [t_2, t_2 + c_d n^{1/d}]$
is at
most $T(d+1, n) = (n^2/2) \left(1 - 1/d \right) + O_d (1)$.
Hence, the number of pairs for which
$d(p_{i(1)}, p_{i(2)})$
belongs to the union of the two intervals in question is at most
$(n^2/2)( 1 - 1/d + 2 \delta ) + O_d (1) $.
In view of \thetag{1.3},
the last expression is bounded from above by
$$
T \left( m(d - 1) + 1, n \right) = \frac{n^2}{2} \left(1 - \frac{1}{m(d - 1)}\right) + O_d (1) ,
$$
provided that
$2 \delta < 1/d - 1/(d + 1) \leq 1/d - 1/m(d - 1)$ and $n$ is
sufficiently large.
Thus, in the case
$t_1 \le  C_d n^{1/d}$, Theorem 1 is true.

(3)
It follows from part (2) that
$$
{\frac{t_2 + c_d n^{1/d}}{t_2}} \le {\frac{t_1 + c_d n^{1/d}}{t_1}} \le
{\frac{C_d n^{1/d} + c_d n^{1/d}}{C_d n^{1/d}}} = 1 + {\frac{c_d}{C_d}} \le 1+c_d,
$$
which proves (3).
\hfill
$\square $


\medskip

{\it{Proof of Theorem 1, continuation.}}
{\bf{Step 2}}.
In the rest of the proof, we assume that conditions
(1), (2), and (3) of Lemma 1 are satisfied.
Consider the {\it{graph $G$ with vertex set $\{p_1, \ldots ,
p_n\}$, where $p_{i(1)}$ and $p_{i(2)}$ are connected by an edge if and only if
$d(p_{i(1)}, p_{i(2)})$ belongs to one
of the two intervals in question.}} 

Suppose, in order to obtain a contradiction, that the number of
edges of $G$ is greater than $T \left( m(d - 1) + 1, n \right)
$.
By \cite
{1}, 
 Ch.\ 6, $G$ contains a subgraph
$G_1 = K(1,1, \ldots,
1, \lfloor \text{\rm const}_d \cdot n\rfloor)$, that is,
a complete $ \left( m(d - 1) + 1 \right) $-partite graph
with $1,1,\ldots, 1, \lfloor \text{\rm const}_d \cdot n
\rfloor $
points in its parts called {\it{primary colour classes}}.
(So here we consider the vertices coloured.) Obviously, it
makes sense to speak about the $j$-th primary
colour class of any (spanned)
subgraph
of $G_1$, for $1 \le j \le m(d - 1) + 1$.
This is
the
intersection of the $j$-th primary colour class of $G_1$
with
the vertex set of the subgraph.
For the subgraphs considered later in this proof, these
primary colour classes are always non-empty.

Define the {\it{secondary colouring}} of the {\it{edges}} of
$G_1$, as follows. Assign to each edge the symbol $L$ and $R$,
according to whether the
length of the corresponding segment lies in the first or in the second
interval.
Since the two intervals are disjoint (cf. Lemma 1 (1)),
the secondary colouring is uniquely determined.

At least half of the points of the $ \left( m(d-1)+1 \right)
$-st primary
colour class of $G_1 = K(1,1, \ldots, 1,
\lfloor \text{\rm const}_d \cdot n\rfloor)$ are
joined by edges of the same secondary
colour $L$ or $R$ to the unique point in the first primary
colour class. By induction, the unique points
in the 1st, 2nd, $\ldots $, $m(d-1)$-st primary colour classes of
$G_1$
and
some $\lfloor \text{\rm const}_d \cdot n\rfloor $ points
in the $ \left( m(d-1)+1 \right) $-st primary
colour class of $G_1$ satisfy
the following.
The secondary colour
of an edge between any two of these points only
depends on the primary colour classes the endpoints of the edge
belong to. We
denote the subgraph induced by all these points by
$G_1^*$.

From now on, we will consider $G_1^*$ rather than $G_1$.
We will show that such a graph $G_1^*$
cannot exist,
for $c_d > 0$
a sufficiently small constant.
This contradiction will prove that
the number of pairs $\{p_{i(1)}, p_{i(2)} \}$,
whose distances lie in the union of our
two intervals, is at most as large as is stated in Theorem 1.

\medskip

{\bf{Step 3}}.
Let $D_d > 2$ be a sufficiently large constant. We distinguish two cases:
$$
{\text{Case I: }}\;\;\;
t_2/(t_1 + c_d n^{1/d}) \le D_d,\;\;\;\;\;\;\;\;\;
$$
\vskip-0.7cm
$$
{\text{Case II: }}\;\; t_2/(t_1 + c_d n^{1/d}) > D_d >
2\,.\;\;
$$

In Case I, we use two-distance sets in ${\Bbb{R}}^{d-1}$.
The proof will be presented in {\bf{Step 4}}, and will be
completed using Lemma 2.

In Case II, the two types of distances, i.e., those belonging to the
first interval and to the second
one, can be treated separately. The segments corresponding to
different types
of distances will turn out to be ``almost orthogonal''.
We will describe the structure of
our edge coloured graph.
The proof in this case will be carried out in {\bf{Step 5}} and
completed by Lemma 9.


\medskip

{\bf{Step 4}}. First, we analyze Case I. The proof of the following
lemma consists of six easy parts (enumerated as {\bf{A}}, {\bf{B}},
$\ldots $, {\bf{F}}).

\medskip


{\bf{Lemma 2.}}
{\it{The upper estimate of Theorem 1 holds in Case I.}}

\medskip

{\it{Proof.}}
{\bf{A.}}
In Case I, we have by Lemma 1 (3), for $c_d > 0$
sufficiently small
$$
\frac{t_2 + c_d n^{1/d}}{t_1} = \frac{t_2 + c_d n^{1/d}}{t_2}
\cdot
\frac{t_2}{t_1 + c_d n^{1/d}} \cdot \frac{t_1 + c_d n^{1/d}}{t_1}
\leq (1+c_d)^2 D_d \le \text{\rm const}_d \cdot D_d \,.
$$
Thus, in Case I,
the quotient of any two distances lying in the union
of our two intervals is at most ${\text{\rm const}}_d \cdot D_d$. Therefore,
these quotients lie between two positive bounds, namely
$(\text{\rm const}_d \cdot
D_d)^{-1}$ and $\text{\rm const}_d \cdot D_d$.
In particular, this holds for the distances between the endpoints of the
edges of the graph $G_1^*$.


\medskip

{\bf{B.}}
\;{\bf{Definition 3.}}
Let $m > d$ be any integer, and let $x_1, \ldots ,
x_m,$ $x_{m+1} \in {\Bbb{R}}^d$ be any distinct points in
${\Bbb{R}}^d$.
Let
$\Delta (x_1, \ldots, x_m, x_{m + 1})$ be the maximum
absolute value of all determinants whose columns are any $d$
vectors from the set
$$
\{(x_1 - x_{m+1})/
d(x_1,x_{m+1}), \ldots , (x_m - x_{m+1})/d(x_m,x_{m+1})\}
\subset S^{d - 1}.
$$


\medskip

Clearly, $\Delta (x_1, \ldots, x_{m+1})$
is invariant under simultaneous similarity
transformations of $x_1, \ldots, x_{m+1}$. Furthermore,
it is nonnegative, and equals $0$
if and only if $x_1, \ldots, x_{m+1}$ lie in an (affine) hyperplane.
Thus, it can be considered as a measure of ``non-hyperplanarity of
$x_1, \ldots, x_{m + 1}$''.
We will apply the above definition
for the case $m := m(d-1)$ (recall \thetag{1.3}).

\medskip

{\bf{C.}}
{\bf{Claim 1.}} {\it{Let $ q_1,\ldots , q_{m(d-1) + 1}$ be
vertices
of $G_1^*$, one from each respective primary colour class.}}
(Thus, $ q_1,\ldots , q_{m(d-1)}$ are fixed, but
$q_{m(d-1) + 1}$ can assume $\lfloor {\text{const}}_d \cdot n
\rfloor $ values, i.e., points.)
{\it{Then
$\Delta(q_1, \ldots, q_{m(d - 1) + 1})$ is at least some
positive constant, provided that $c_d > 0$
is sufficiently small.}}


\medskip

{\it{Proof.}} Suppose, for contradiction, that $c_d>0$ is very small,
i.e., we
have $c_d^N < 1/N$, say, for a large integer $N$, but
$\Delta(q_1,\ldots,
q_{m(d - 1) + 1})
$ can get arbitrarily close to $0$.
That is, $\Delta(q_1, \ldots,
q_{m(d - 1) + 1}) <
1/N$, say, for some choice $q_1^N, \ldots , q_{m(d - 1) + 1}^N$ of the
points $q_1, \ldots, q_{m(d - 1) + 1}$. (Actually, only the
last point can vary.)
We apply to each of $q_1^N, \ldots, q_{m(d - 1) + 1}^N$,
simultaneously,
a similarity transformation $\Phi _{\lambda }$ with ratio $\lambda  >  0$
such that the following holds. We have
${\text{\rm{diam}}}\,\{\Phi _{\lambda }q_1^N, \ldots, \Phi _{\lambda }
q_{m(d - 1) + 1}^N\}
= 1$ and $\{ \Phi _{\lambda }q_1^N, \ldots, \Phi _{\lambda }
q_{m(d - 1)
+ 1}^N \} $ lies in the unit ball of ${\Bbb{R}}^d$. Then, by {\bf{A}},
the minimal distance in
$\{ \Phi _{\lambda } q_1^N, \ldots, \Phi _{\lambda }
q_{m(d - 1) + 1}^N \} $
is at least
$1/({\text{\rm const}}_d \cdot D_d) > 0$.

Now let $N \to \infty$.
Then, for a certain subsequence $N(\nu )$ of the $N$'s,
the following four statements are true:
\newline
\noindent
(i)~for each $1 \le j \le m(d-1) + 1$, we have that
$\lim _{\nu \to \infty }\Phi _{\lambda }(q_j^{N(\nu )})$ exists;
\newline
\noindent
(ii)~these limit points have pairwise distances at least
$1/({\text{\rm const}}_d \cdot D_d)$;
\newline
\noindent
(iii) for any $j(1) \neq j(2)$, the distance
$d(\Phi _{\lambda }(q_{j(1)}^{N(\nu )}), \Phi _{\lambda }
(q_{j(2)}^{N(\nu )}))$ lies
in $[ \lambda t_{\kappa }, \lambda (t_{\kappa } + c_d n^{1/d})]$
for some $\kappa \in \{ 1,2 \} $.
\newline
\noindent
(iv) $\lim _{\nu \to \infty } \Delta \left(
q_1^{N(\nu )}, \ldots,
q_{m(d - 1) + 1}^{N(\nu )} \right) = 0$.

By (ii),
$\left( \lim _{\nu \to \infty }\Phi _{\lambda }
(q_1^{N(\nu )}), \ldots ,
\lim _{\nu \to \infty }\Phi _{\lambda }
(q_{m(d-1)+1}^{N(\nu )}) \right)
$
belongs to the domain of
definition of the function $\Delta (x_1, \ldots ,
x_{m(d - 1) + 1})$.

By (iii), Lemma 1 (3)
and $c_d^{N(\nu )} < 1/N(\nu )$,
we have that
any two numbers that belong to the same new interval
$[ \lambda t_{\kappa }, \lambda (t_{\kappa } + c_d n^{1/d}) ]$,
for
$\kappa  \in \{ 1,2 \}$,
have a ratio in
$[(1 + c_d^{N(\nu )}
)^{-1}, 1 + c_d^{N(\nu )}
] \subset
[\left( 1 + (1/N(\nu ))
\right) ^{-1}, 1 + (1/N(\nu ))
]$. 
Thus, this ratio lies
in an as small neighbourhood of $1$, as we want.
Therefore, for $\nu \to \infty $, both our
new $\kappa $-th
intervals
converge to degenerate intervals, i.e.,\ to points. In particular, the
second new interval converges to $\{ 1 \} $.

By the similarity invariance of $\Delta( \cdot )$ and (iv),
we have
$$
\Delta \left( \lim _{\nu \to \infty }\Phi _{\lambda }
q_1^{N(\nu )}, \ldots,
\lim _{\nu \to \infty }\Phi _{\lambda }q_{m(d - 1) + 1}^
{N(\nu )} \right) 
$$
$$
=
\lim _{\nu \to \infty } \Delta \left( \Phi _{\lambda }
q_1^{N(\nu )}, \ldots,
\Phi _{\lambda }q_{m(d - 1) + 1}^{N(\nu )} \right) = 0.
$$
That is, the points $\lim _{\nu \to \infty }\Phi _{\lambda }q_1^{N(\nu )},
\ldots, \lim _{\nu \to \infty }
\Phi _{\lambda }q_{m(d - 1) + 1}^{N(\nu )}$ lie in some
hyperplane of ${\Bbb{R}}^d$,
their number is $m(d-1) + 1$, and they determine
only two distinct distances. This contradiction ends the proof of Claim 1.
\hfill
$\square $


\medskip

{\bf{D.}}
Let us fix some $q_{m(d-1)+1}$ in the $\left( m(d-1)+1 \right)
$-st primary colour
class of $G_1^*$.
By Claim 1, among the $m(d-1) > d$ unit vectors
$$
u_1(q_{m(d-1)+1}) :=
(q_1 - q_{m(d-1)+1}) / d(q_1,q_{m(d-1)+1}),
\ldots ,
$$
$$
u_{m(d-1)}(q_{m(d-1)+1}) :=
(q_{m(d-1)} - q_{m(d-1)+1}) / d(q_1,q_{m(d-1)+1}),
$$
there are
$u_{j(1)}(q_{m(d-1)+1}), \ldots , u_{j(d)}(q_{m(d-1)+1})$
such that
$$
| {\text{det}}
\left( u_{j(1)}(q_{m(d-1)+1}),
\ldots , u_{j(d)}(q_{m(d-1)+1}) \right) | \ge
{\text{const}}_d > 0 .
$$
Since there are only ${\text{const}}_d$ choices for these
$d$-tuples, still for $\lfloor {\text{const}} _d \cdot n
\rfloor $ many
choices of $q_{m(d-1)+1}$ this $d$-tuple is the same,
$\{ u_{j(1)}, \ldots , u_{j(d)} \} $, say. We
will write $C_{m(d - 1) + 1}$
for the set of these $\lfloor {\text{const}}_d \cdot n
\rfloor $
points $q_{m(d-1)+1}$
in the $\left( m(d-1)+1 \right) $-st primary colour
class of $G_1^*$.
Thus,
$$
| {\text{det}}  \left( u_{j(1)}, \ldots , u_{j(d)} \right)|
\ge
{\text{const}}_d > 0\,.
$$
Hence, we have a {\it{$d$-dimensional simplex $S(q_{j(1)},
\ldots , q_{j(d)}, q_{m(d-1)+1})$}},
and a {\it{$(d-1)$-dimensional simplex $S(q_{j(1)},
\ldots , q_{j(d)})$}}.

From the set of unit vectors $u_1, \ldots , u_{m(d-1)}$,
we will consider only $u_{j(1)}, \ldots ,$
\newline
$ u_{j(d)}$.
Further, from among all ${\text{\rm{const}}} _d \cdot n$
points $q_{m(d-1)+1}$ in the $ \left( m(d-1)+ 
1 \right) $-st primary colour class of $G_1^*$,
we will restrict our attention to the subset $C_{m(d - 1) + 1}$.
{\it{We write $G_1^{*}{'}$ for the induced
subgraph of $G_1^*$, containing all
(single) vertices of $G_1^*$ in its first $m(d - 1)$ primary colour classes,
and $C_{m(d - 1) + 1}$ from its last primary colour class.}}


\medskip
{\bf{E.}}
Recall from {\bf{Step 2}} of the proof of Theorem 1 the
following.
For any $h = 1, \ldots, d$, either
\newline
(1) for any choice of the
vertex $q_{m(d - 1) + 1,i}$ of the $ \left( m(d - 1) + 1
\right) $-st
primary colour class of
$G_1^*$, the distance
$d(q_{j(h)}, q_{m(d - 1) + 1,i}) $
lies in the first
interval, or
\newline
(2) for any choice of the vertex
$q_{m(d - 1) + 1,i}$ of
the $ \left( m(d - 1) + 1 \right) $-st primary colour class of
$G_1^*$, the distance $d(q_{j(h)}, q_{m(d - 1) + 1,i}) $
lies in the second interval. In particular, this holds for the
subgraph $G_1^{*}{'}$ of $G_1^*$.
This means that
$$
C_{m(d - 1) + 1} \subset \cap _{h = 1} ^d S_{j(h)},
\tag 3.1
$$
where
$S_{j(1)}, \ldots ,S_{j(d)}$ are
spherical shells
with centres $q_{j(1)},
\ldots, q_{j(d)}$, inner radii either
$t_1$ or $t_2$, and outer radii either
$t_1 + c_d n^{1/d}$ or $t_2 + c_d n^{1/d}$, resp. For
each of these spherical shells, the quotient of the difference of the outer
and inner radii and of the inner radius is
$c_dn^{1/d}/t_{\kappa } \le c_dn^{1/d}/
(C_dn^{1/d}) = c_d/C_d \le c_d$, by Lemma 1
(2). Hence, this quotient is
in an arbitrarily small
neighbourhood of $0$, if $c_d > 0$
is chosen sufficiently small.

We are going to show that, for $c_d > 0$
sufficiently small, the inclusion
\thetag{3.1}
is impossible, yielding
the desired contradiction.

Before this, we have to introduce some notations.
Observe that ${\text{aff}}\, \{ q_{j(1)}, \ldots, $
\newline
$q_{j(d)} \} $
is a hyperplane of symmetry of $\cap _{h = 1}^d S_{j(h)}$, which
will be identified with the hyperplane $x_d=0$. (As
$S(q_{j(1)}, \ldots , q_{j(d)})$
is $(d-1)$-dimensional, so is its affine hull.)
Let $H^+$ and $H^-$ denote the closed half-spaces $x_d \ge 0$
and $x_d \le 0$, resp. One of them contains at least
half of the points of
$C_{m(d - 1) + 1}$. We may suppose this is
$H^+$. Thus
$$
H^+ {\text{ contains }}
\lfloor {\text{const}}_d \cdot n \rfloor {\text{ points of }}
C_{m(d - 1) + 1}.
$$
Let us fix a point 
$q_{{m(d - 1) + 1},1} \in C_{m(d - 1) + 1} \cap H^+$.

For any $h \in \{1, \ldots , d \} $,
define the slab $S'_{j(h)}$, as follows. Let
$S'_{j(h)}$ be bounded by two hyperplanes, both orthogonal to
$q_{j(h)} - q_{m(d-1)+1,1}$.
Further, they intersect the half-line from $q_{j(h)}$, passing through
$q_{m(d-1)+1,1}$, at points with distances $d(q_{j(h)},
q_{m(d-1)+1,1}) + c_d
n^{1/d}$ and $d(q_{j(h)},q_{m(d-1)+1,1}) - 2 c_dn^{1/d}$ from
$q_{j(h)}$. (By Lemma 1 (2), this difference is positive, for
$c_d > 0$ sufficiently small.)
We need the following

\medskip


{\bf{Claim 2.}} {\it{If $c_d > 0$ is sufficiently small, then
$$
(\cap _{h = 1}^d S_{j(h)})\cap H^+
\subset \Pi := \cap _{h = 1}^d S'_{j(h)} .
$$
This holds both in Case I and in Case II.}}


\medskip

This statement appears to be intuitively clear, but we have been
unable to show it by a simple geometric argument. We provide a
proof in the original version
of our paper \cite
{10}, 
on arXiv;
see parts {\bf{8}}-{\bf{10}}
of the proof of Theorem 1,
pp. 9-13. It uses elements of the algebraic topology of
Euclidean spaces \cite
{17}. 


\medskip

{\bf{F.}} Again, we handle both Cases I and II. 
The set $\Pi $ in Claim 2 is a parallelepiped,
and is circumscribed about a ball of
diameter $3c_d n^{1/d}$.
Its volume is
$(3c_d/2)^d n$ times the volume of its homothetic copy $\Pi ^1$
circumscribed about the unit ball. Moreover,
$$
V(\Pi ^1) =
2^d/|{\text{det}} \,(u_{j(1),1}, \ldots , u_{j(d),1})|,
$$
with the denominator at least
${\text{const}}_d > 0$, by {\bf{D}}.
(The easiest way to see this volume formula
is as follows. The polar body
$(\Pi ^1)^*$ of
$\Pi ^1$ is a cross-polytope, with $V\left( (\Pi ^1)^*
\right) = (2^d/d!)|{\text{det}} \,(u_{j(1),1},
\ldots , u_{j(d),1})|$.
Simultaneously, the product $V(\Pi ^1)V \left(
(\Pi ^1)^* \right) $
of the two
volumes is invariant under linear maps. Hence, it equals
$4^d/d!$, as can be
calculated from the case when $\Pi _1$ is the unit cube.
Cf. \cite
{20}, 
 pp. 165, 169.)
All these imply that
$$
V(\Pi ) \le (3c_d/2)^d n 2^d / {\text{const}}_d \,.
$$

Now a standard volume consideration finishes the proof of Lemma 2.
Consider the open
balls of unit diameter, with centres at all
$\lfloor {\text{const}}_d \cdot n  \rfloor $
points $q_{m(d-1)+1,i} \in C_{m(d - 1) + 1} \cap H^+
\subset \Pi $.
{\it{These are pairwise disjoint open balls contained in a
concentric homothetic copy $\Pi '$
of $\Pi $, with inradius $3c_dn^{1/d}/2 + 1/2$.}}
However,
$$
V(\Pi ') \le {\text{const}}_d \cdot c_d ^d n .
$$
So, if $c_d>0$ is sufficiently small, then the volume of $\Pi'$ is not
large enough to contain 
$\lfloor {\text{const}}_d \cdot n \rfloor $ disjoint open
balls of unit diameter. This contradiction completes the proof
of Lemma 2 and, hence,
Theorem 1 in Case I (see {\bf{Step 3}}). \hfill
$ \square $

\medskip

Next, we turn to the proof of Theorem 1 in Case II. 

\medskip

{\it{Proof of Theorem 1, continuation.}} {\bf{Step 5}}. Now we
assume that $t_2/(t_1 + c_d n^{1/d}) > D_d > 2$, where $D_d$ is a
sufficiently large constant (Case II).

We investigate the secondary (edge)
colourings of the graph $G_1^{*}{'}$
from {\bf{Step 2}} of the proof of Lemma 2. Each edge is coloured
either by $L$ or by $R$.
Each edge coloured by $R$ has length at least $t_2$,
and each edge coloured by
$L$ has length at most $t_1 + c_d n^{1/d}$.
By $t_2/(t_1 + c_d n^{1/d}) > 2$,
{\it{any edge
coloured by $R$ is more than twice as long as any edge
coloured by $L$}}.

This implies that one can define an equivalence relation $\sim $
on the vertices of $G_1^{*}{'}$ as follows.

\medskip


{\bf{Definition 4.}} 
For any two vertices $q_{j(1)}$, $q_{j(2)}$ of $G_1^{*}{'}$,
we write $q_{j(1)} \sim q_{j(2)}$ if either
$q_{j(1)} = q_{j(2)}$, or the edge $q_{j(1)} q_{j(2)}$ is
coloured by $L$. 

\medskip


Recall from {\bf{Step 2}} and the proof of Lemma 2, {\bf{D}}, 
that the colour of an edge of $G_1^{*}{'}$
between vertices of two primary colour classes does not depend
on the vertices chosen from the primary colour classes. (This
is equivalent to its special case when one of the
primary colour classes is the $\left( m(d-1) + 1 \right )$-st
primary
colour class.)
$$
{\text{Therefore, we may consider the relation }} \sim
{\text{ as defined}}
$$
$$
{\text{alternatively on the set of primary colour classes of }}
G_1^{*}{'}.
$$
Whether we consider it on the set of vertices, or on the
primary colour classes, will be
clear from the context.
Let $\ell $ denote the number of $\sim $-equivalence classes.

Let us choose for each of the $\ell  
\sim $-equivalence classes of primary
colour classes of $G_1^{*}{'}$ one vertex $q_j$
from their union; let these be $r_1, r_2, \ldots, r_{\ell }$.
By Lemma 1 (3),
edges coloured the same way have a ratio in an as
small neighbourhood of~1 as we want, provided $c_d > 0$
is sufficiently small.
Note that any edge among $r_1, \ldots, r_{\ell }$
is coloured by $R$. Therefore the quotients of
the lengths of these edges are in an as small neighbourhood of $1$ as we want,
for $c_d > 0$ sufficiently small.
Hence, for $c_d > 0$ sufficiently small,
we have $\ell \leq d + 1$. Namely, for $d
+ 2$ points in ${\Bbb{R}}^d$
the quotient of the maximum and the minimum distances is
at least some constant strictly greater than~$1$. (Cf.\
Sch\"utte \cite
{24}, 
 Satz 3, which gives the sharp lower bound, which is
$(1 + 2/d) ^{1/2}$, for $d$ even and
$[1 + 2(d + 2)/\left( d(d + 2) - 1 \right) ] ^{1/2}$, for $d$
odd.)
The same argument shows that $r_1, \ldots , r_{\ell }$ cannot lie
in an affine $(\ell - 2)$-plane, thus determine an
$(\ell - 1)$-simplex, namely $S(r_1, \ldots , r_{\ell })$.


Our goal is
to show that the simplex
$S(r_1, \ldots ,r_{\ell })$ is ``close''
to a regular
$(\ell - 1)$-simplex of edge length $t_2$. Similarly,
the vertices of
$G_1^{*}{'}$
in single
$\sim $-equivalence classes
are ``close'' to the vertices of regular simplices of
edge length $t_1$, of dimensions at most $d - \ell + 1$, with affine
hulls nearly
orthogonal to ${\text{aff}}\, \{ r_1, \ldots , r_{\ell } \} $.
The number of primary colour classes of $G_1^{*}{'}$ is maximum
if all of the last simplices have dimension $d - \ell + 1$.


\medskip
{\bf{Lemma 3.}}
{\it{In Case II, the number $\ell $ of the
$\sim$-equivalence classes is at least
$2$, provided $c_d > 0$ is sufficiently small.}}

\medskip

{\it{Proof.}}
If $\ell = 1$, then all distances between the vertices of
$G_1^{*}{'}$
lie in $[t_1, t_1 + c_d n^{1/d}]$, contradicting Lemma 1 (1).
This proves Lemma 3.
\hfill
$\square $


\medskip

{\bf{Lemma 4.}}
{\it{In Case II, let $q_j \in V(G_1^{*}{'})$ be
in the $\sim $-equivalence class of $r_1 \in V(G_1^{*}{'})$
such that $q_j \neq
r_1$. Further, let $r_2 \in V(G_1^{*}{'})$ be in another
$\sim $-equivalence class, as $r_1$.
Then
$
|\sphericalangle q_j r_1 r_2 - \pi /2|
$
is as small as we want, for
$D_d$ sufficiently large and $c_d > 0$ sufficiently small.}}
(Here $r_2$
exists by Lemma 3.)

\medskip


{\it{Proof.}}
We are
going to estimate from above
$$
|\cos(\sphericalangle q_j r_1 r_2)|
= \left| d(r_1, r_2)^2
+ d(r_1, q_j)^2 - d(q_j, r_2)^2 \right|
/
\left( 2 d(r_1, r_2) d(r_1, q_j) \right)
$$
$$
=\left| \left( d(r_1, r_2) + d(q_j, r_2) \right) \left(
d(r_1, r_2)
- d(q_j, r_2) \right) + d(r_1, q_j)^2 \right|
$$
$$
/ \left( 2 d(r_1, r_2) d(r_1, q_j) \right)
\leq [2(t_2 + c_d n^{1/d}) c_d n^{1/d} + (t_1 + c_d n^{1/d})^2]
/(2 t_1 t_2)\,.
$$
By Lemma 1 (3), any two numbers from the same
interval $[t_{\kappa }, t_{\kappa } + c_d n^{1/d}]$
have quotients as close to $1$ as we want,
for $c_d > 0 $
sufficiently small.
Thus, we suppose $t_1 + c_d n^{1/d} \leq 2t_1$ and
$t_2 + c_d n^{1/d} \leq 2 t_2$,
for $c_d > 0$ sufficiently small.
Then
$$
|\cos(\sphericalangle q_j r_1 r_2)|
\leq [2 \cdot 2 t_2 \cdot c_d n^{1/d} + 4 t_1^2]
/(2 t_1 t_2) =
$$
$$
2 c_d n^{1/d} / t_1 + 2 t_1 / t_2 \leq 2 c_d/ C_d + 2/ D_d
< 2 c_d + 2/ D_d ,
$$
by Lemma 1 (2), and by $t_2 / t_1 \ge t_2 / (t_1 + c_d n^{1/d}) > D_d$
(Case II).
If $D_d$
is sufficiently large and $c_d > 0$
is sufficiently small,
then this last expression, and hence also
$|\cos(\sphericalangle q_j r_1 r_2)|$ is as small as we want.
This proves Lemma 4.
\hfill
$\square $


\medskip

{\bf{Lemma 5.}}
{\it{In Case II, the number $\ell $ of the $\sim $-equivalence
classes is at most
$d$, for $D_d$ sufficiently
large and $c_d > 0$
sufficiently small.}}

\medskip


{\it{Proof.}}
We already know that
$\ell \leq d + 1$ (cf.\ {\bf{Step 5}}),
so we have to exclude $\ell = d + 1$ only.

Suppose $\ell = d + 1$. At the beginning of {\bf{Step 5}}, we
selected
points $r_1, \ldots, r_{\ell } = r_{d + 1}$, one from the union of each
$\sim $-equivalence class of the primary colour classes of
$V(G_1^{*}{'})$.
By Lemma 1 (3),
for $c_d > 0$
sufficiently small, we have that the quotients of any two distances
among these points are in an as small neighbourhood of $1$ as we want.
This
also implies that any three of these points determine a triangle with angles
as close to $\pi /3$ as we want.

Let $v_{\mu } := (r_{\mu } - r_{d + 1})  / d(r_{\mu }, r_{d + 1}) \in S^{d - 1}$,
for $\mu = 1, \ldots, d$.
Let $V$ denote the $d \times d$ matrix with columns
$v_1, \ldots, v_d$.
We have
$$
|{\text{det}}\,(v_1, \ldots, v_d)| = |{\text{det}}\, V|
= [ {\text{det}}\,(V'V) ] ^{1/2}
= [ {\text{det}}\,(\langle v_{\mu (1)}, v_{\mu (2)}\rangle) ]
^{1/2} ,
$$
where $V'$ is the transposed matrix of $V$. Moreover,
$(\langle v_{\mu (1)}, v_{\mu (2)}\rangle)$ is a $d \times
d$ matrix for which
$\langle v_{\mu (1)}, v_{\mu (1)} \rangle = 1$, and
$\mu (1) \neq \mu (2)$ implies that
$\langle v_{\mu (1)}, v_{\mu (2)} \rangle $ is as
close to $\cos (\pi /3) = 1/2$ as we want.
Hence, $|{\text{det}}\,(v_1, \ldots, v_d)|$ is as close
to $[{\text{det}} \left( (1 + \right.
$
\newline
$
\left. \delta_{\mu (1) \mu (2)})/2 \right) ]^{1/2}$
as we want,
for $D_d$ sufficiently large and $c_d > 0$
sufficiently small.
Here $[{\text{det}}\left( (1 + \delta_{\mu (1) \mu (2)})/2 \right)]
^{1/2}$ equals
the absolute
value of the determinant whose co\-lumns are the
unit vectors pointing from a vertex of a
regular $d$-simplex to all other vertices. Thus it is a non-zero
constant.
All this implies that
$|{\text{det}}\,(v_1, \ldots, v_d)|$ {\it{is greater than a
non-zero constant}}. In particular,
$v_1,\ldots, v_d$ {\it{are linearly independent}}.

By \thetag{1.3}, $m(d - 1) + 1 > d + 1$ for $d \ge 2$,
hence some of the $\ell = d + 1$ 
$\sim $-equivalence classes
must contain at least two
vertices $q_j$ of $G_1^{*}{'}$. Assume without loss of
generality that $r_{d+1}$ belongs to such a class and $q_j$ is
one of its elements different from $r_{d+1}$.
In view of Lemma 4, the scalar product of the vector $v:=
(q_j - r_{d + 1}) / d(q_j, r_{d + 1}) \in S^{d - 1}$
with any $v_{\mu}$, $1 \le \mu \le d$, is as close to $0$ as we
want, provided that $D_d$ is sufficiently large and
$c_d > 0$ is sufficiently small.
In other words, $\max _{1 \le \mu \le d} |\langle v, v_{\mu }
\rangle |$
{\it{is as
small as we want, for $D_d$ sufficiently large and
$c_d > 0$ sufficiently small}}.

By the linear independence of
$v_1, \ldots, v_d$, we have $v = \sum _{\mu = 1}^d
\lambda_{\mu } v_{\mu }$ for some
$\lambda_ {\mu } \in {\Bbb{R}}$.
Consider this as a system of equations for $\lambda_ {\mu }$,
and note that the absolute value of any coordinate of $v$ and any
$v_{\mu}$ is at most $1$.
Then we have by Cramer's rule for
$\lambda_{\mu }$, and by
$|{\text{det}}\,(v_1, \ldots, v_d)| \ge {\text{const}}_d$,
that $|\lambda_{\mu }| \leq {\text{\rm const}}_d/
|{\text{det}}\,(v_{\mu })|
\leq {\text{\rm const}}_d$.

Hence
$$
1 = \langle v, v \rangle
= \langle v , \sum _{\mu = 1}^d \lambda_{\mu } v_{\mu } \rangle
= \sum _{\mu = 1}^d \lambda_{\mu } \langle v, v_{\mu } \rangle
\leq
$$
$$
d \cdot \max _{1 \le \mu \le d}|\lambda_{\mu }| \cdot
\max _{1 \le \mu \le d}
|\langle v, v_{\mu }\rangle| \leq 
\text{\rm const}_d
\cdot \max _{1 \le \mu \le d}
| \langle v, v_{\mu } \rangle | .
$$
This contradicts the fact that $\max _{1 \le  \mu \le d}
|\langle v, v_{\mu } \rangle |$ is as small as we want,
for $D_d$ sufficiently large and
$c_d > 0$ sufficiently small. This completes the proof of Lemma
5. 
\hfill
$\square $

\medskip


{\bf{Lemma 6.}}
{\it{In Case II, for $D_d$ sufficiently large and
$c_d > 0$ sufficiently small,
any $\sim $-equivalence class contains at most $d - \ell + 2$
points~$q_j$.}}

\medskip


{\it{Proof.}}
Let $q_j$ be any vertex in the $\sim $-equivalence
class of $r_\ell $, say, and
let
$q_j \neq r_\ell $. (If such a vertex did not exist, then
this $\sim $-equivalence class
would have $1 < d - \ell + 2$ points, by Lemma 5.)
Then, by Lemma 4,
$|\cos(\sphericalangle q_j r_\ell r_1)|, \ldots,
|\cos(\sphericalangle q_j r_\ell r_{\ell - 1})|$
are as small as we want, for $D_d$ sufficiently large and
$c_d > 0$ sufficiently small.
Let $w^j := (q_j - r_\ell )/ d(q_j, r_{\ell }) \in S^{d - 1}$,
and, for $\mu = 1, \ldots, \ell - 1$, let $w_{\mu }
:= (r_{\mu } - r_{\ell })/ d(r_{\mu }, r_\ell ) \in S^{d - 1}
$. Then we have that
$| \langle w^j, w_{\mu } \rangle |$
{\it{is as small as we want, for
$D_d$ sufficiently large and $c_d > 0$ sufficiently
small}}.
{\it{Suppose $r_\ell = 0$}}, and let
$(w^j)'$ be the orthogonal projection of $w^j$ to
the linear $(\ell - 1)$-subspace
$\text{\rm
aff} \, \{r_1, \ldots, r_\ell \} = \text{\rm
lin} \, \{r_1, \ldots, r_ {\ell - 1} \} $.
Then
$w_1, \ldots, w_{\ell -
1}, (w^j)' \in \text{\rm lin} \, \{r_1, \ldots,
r_{\ell - 1}\}$. Moreover,
{\it{$| \langle (w^j)', w_{\mu } \rangle |
= | \langle w^j,  w_{\mu } \rangle |$ is as small as we
want, for $D_d$ sufficiently large and
$c_d > 0$ sufficiently small.}}

Now we proceed in the linear subspace $\text{\rm lin} \,
\{r_1, \ldots, r_{\ell - 1} \} $,
as we proceeded in ${\Bbb{R}}^d$ in the proof of Lemma 5.
We have $|{\text{det}}\,(w_1, \ldots, w_{\ell -
1})| \geq {\text{const}}_d  > 0$.
(Observe that $\ell $ can assume
only finitely many values. This is why we could write here
${\text{const}}_d > 0$.)
Moreover, $(w^j)' =
\sum _{\mu = 1}^{\ell -1} \lambda _{\mu } w_{\mu }$,
where now by Cramer's rule $| \lambda _{\mu } | \le {\text{const}}_d \cdot \| (w^j)' \| $. Then
$$
\| (w^j)' \| ^2 = \langle (w^j)', (w^j)' \rangle =
\langle  (w^j)', \sum _{\mu =1}^{\ell - 1}
\lambda _{\mu }w_{\mu } \rangle =
\sum _{\mu =1}^{\ell - 1}
\lambda_{\mu } \langle (w^j)', w_{\mu } \rangle 
$$
$$
\leq d \cdot \max _{1 \le \mu \le \ell - 1}
|\lambda_{\mu }| \cdot
\max _{1 \le \mu \le \ell - 1}
|\langle (w^j)', w_{\mu } \rangle|
$$
$$
\leq \text{\rm const}_d \cdot \| (w^j)' \|\cdot
\max _{1 \le \mu \le \ell - 1}
|\langle  (w^j)', w_{\mu } \rangle | .
$$
Hence,
$$
\| (w^j)' \| \leq \eta := \text{\rm
const}_d \cdot \max _{1 \le \mu \le \ell - 1}
|\langle  (w^j)', w_{\mu } \rangle |,
$$
and here
$\eta $
{\it{is as small as we want, for $D_d$
sufficiently large and $c_d > 0$ sufficiently small}}.

Suppose that the equivalence class of $r_\ell $ contains
$d - \ell + 2$ other points, 
$q_{\ell + 1}, \ldots, $
\newline
$q_{d + 2}$, besides $r_{\ell }$
(any of which could be identical
with the point denoted by $q_j$ at the beginning of the proof
of the lemma).
Let $q_{\ell } = q_{\ell }^* := r_\ell = 0$. Further, let
$q_{\ell + 1}^* ,\ldots, q_{d + 2}^*$
denote the orthogonal projections of $q_{\ell + 1}, \ldots, q_{d + 2}$
to the linear $(d - \ell + 1)$-subspace which is the
orthocomplement of \,$\text{\rm lin} \, \{r_1, \ldots, r_{\ell - 1} \} $. Then we have for distinct $j(1),j(2) 
\in \{\ell, \ell + 1, \ldots, d + 2\}$ that
$$
d(q_{j(1)}^*, q_{j(2)}^*) \leq
d(q_{j(1)}, q_{j(2)}) \leq t_1 + c_d n^{1/d}.
$$
On the other hand, for $j(1) \in \{\ell, \ell + 1, \ldots,
d + 2\}$,
the orthogonal
projection of $q_{j(1)}$ to $\text{\rm lin} \,
\{r_1, \ldots, r_{\ell - 1} \} $
is $q_{j(1)} - q_{j(1)}^*$. Here for $j(1) \ge \ell + 1$ we have
$q_{j(1)} = q_{j(1)} - r_\ell = d(q_{j(1)}, r_\ell )
{w^{j(1)}}$, hence
its orthogonal
projection to $\text{\rm lin} \, \{r_1, \ldots,
r_{\ell - 1} \} $ is
$q_{j(1)} - q_{j(1)}^* = d(q_{j(1)}, r_\ell ) (w^{j(1)})'$.
Therefore, we have $d(q_{j(1)} , q_{j(1)}^*) =
d(q_{j(1)}, r_\ell ) \cdot \| (w^{j(1)})' \| \le
d(q_{j(1)}, r_\ell ) \eta $
and, analogously, $d(q_{j(2)} , q_{j(2)}^*) \le
d(q_{j(2)}, r_\ell ) \eta $.
For $j(1) = \ell $, we have $d(q_{j(1)},q_{j(1)}^*) = 0 \le
d(q_{j(1)}, r_{\ell }) \eta = 0$ and, analogously for
$j(2) = \ell $.

These imply by Lemma 1 (2),
for $D_d$ sufficiently large and
$c_d > 0$ sufficiently small, that for $j(1), j(2) \in \{ \ell ,
\ell + 1, \ldots , d+2 \} $ we have (both for $j(\kappa ) =
\ell $, and for $j(\kappa ) > \ell $) that
$$
d(q_{j(1)}^* ,q_{j(2)}^*) \geq d(q_{j(1)}, q_{j(2)}) -
d(q_{j(1)}^*,  q_{j(1)}) - d(q_{j(2)}^*, q_{j(2)}) \geq 
$$
$$
t_1 - \eta \cdot d(q_{j(1)}, r_\ell )
- \eta \cdot d(q_{j(2)}, r_\ell )
\geq t_1 - 2\eta (t_1 + c_d n^{1/d}) \ge t_1 - 4 \eta t_1.
$$

Let $D_d$ be sufficiently large and $c_d > 0$ be
sufficiently small. Then $c_d n^{1/d}/t_1$ is
sufficiently small (Lemma 1 (2)), and also $\eta $ is
sufficiently small.
Therefore, each $d(q_{j(1)}^*, q_{j(2)}^*)$ lies in an interval, whose maximum and
minimum have a quotient that is as close to $1$ as we want.
Thus, there are $d - \ell + 3$ points,
$q_{\ell }^*,q_{\ell +1}^*, \ldots, q_{d + 2}^*$, with this
property in a $(d - \ell + 1)$-dimensional linear subspace of
${\Bbb{R}}^d$;
namely in
the orthocomplement of $\text{\rm lin} \, \{r_1, \ldots, r_{\ell - 1} \}$.
This is
impossible by the theorem of Sch\"utte \cite
{24}, 
 cited at the
beginning of {\bf{Step 5}}, for $D_d$ sufficiently large and
$c_d > 0$
sufficiently small.
\hfill
$\square $


\medskip

{\bf{Lemma 7.}}
{\it{In Case II, for $d \ge 6$, the upper estimate of Theorem 1
holds.}}

\medskip


{\it{Proof.}}
By Lemma 6, the number of all primary colour classes
of $G_1^{*}{'}$, i.e.\ $m(d - 1) + 1$, is at most
$\ell \cdot(d - \ell + 2) \leq \lfloor \left( (d + 2)/2 \right)
^2 \rfloor $.
Hence, by {\thetag{1.1}}, we have $d(d-1)/2 + 1 \leq
m(d - 1) + 1 \leq
\lfloor \left( (d + 2)/2
\right) ^2 \rfloor $.
Thus, $d(d - 1)/2 + 1 \leq \left( (d + 2)/2 \right) ^2$,
implying $d \leq 6$. That is, for $d \ge 7$ we have a contradiction.

For $d = 6$, by \thetag{1.2} we have
$17 = m(d - 1) + 1 \leq \ell(d - \ell +
2) \leq \lfloor \left( (6 + 2)/2 \right)^2 \rfloor = 16$, a
contradiction.

At the
beginning of {\bf{Step 2}}, we assumed, in order to obtain a
contradiction, that Theorem 1 was false. 
This led to a contradiction for every $d \ge 6$.
\hfill
$\square $


\medskip

{\bf{Lemma 8.}}
{\it{In Case II, for $d=2,3$, the upper estimate of Theorem 1 holds.}}


\medskip

{\it{Proof.}}
By \thetag{1.2}, for $d = 2$ we have
$$
4 = m(d - 1) + 1 \le \ell (d - \ell + 2)
\le \lfloor \left( (d + 2)/2 \right) ^2 \rfloor = 4,
$$
implying $\ell = 2$, while for $d = 3$ we have
$$
6 = m(d - 1) + 1
\le \ell (d - \ell + 2) \le
\lfloor \left( (d + 2)/2 \right)^2\rfloor = 6,
$$
implying $\ell \in \{ 2,3 \} $.
For both $d=2$ and $3$, equality in the first inequality implies
that each of the $\ell$ $\sim $-equivalence classes contains
maximally many,
i.e., $d - \ell + 2$ primary colour classes of $G_1^{*}{'}$.

First, let $d = 2$. Then $\ell = d - \ell + 2 = 2$.
Let the $\sim $-equivalence classes on the primary colour classes of
$G_1^{*}{'}$
be
represented by the vertices
$\{ q_1, q_3 \} $ and
$\{ q_2, q_4 \} $.
Here all $q_j$'s belong to
distinct ones among the four primary colour classes
of $G_1^{*}{'}$, which have $1,1,1,\lfloor {\text{const}} \cdot n \rfloor $
vertices, resp.
(Thus, these vertices form a subgraph of
$G_1^{*}{'}$ which is a four-clique --
in particular, each distance determined by them lies in
$[t_1, t_1 + c_d n^{1/d}] \cup [t_2, t_2 + c_d n^{1/d}]$.)
Further,
the secondary
colour of an edge only depends on the primary colour
classes the edge endpoints belong to.
Up to
notation, we may assume that $q_4$ belongs to the last one of these
primary colour classes (it plays the role of $q_{m(d - 1) + 1}$
from part {\bf{E}} of the proof of Lemma 2).
By Lemma 4, $\sphericalangle q_2 q_4 q_3 $ is
close to $ \pi /2 $.
Fix $q_1, q_2, q_3$, and vary $q_4$ in its own
primary colour class in $G_1^{*}{'}$, so that
it assumes $\lfloor {\text{const}} \cdot n \rfloor $
values (points). Then all these points lie
in the intersection of two circular shells
(defined analogously as in part
{\bf{E}} of the proof of Lemma 2).
These have
centres $q_2$ and $q_3$, inner radii some $t_{\kappa }$'s,
and outer radii the respective $(t_{\kappa } +
c_dn^{1/d})$'s.
Moreover,
the unit vectors pointing from $q_4$ to $q_2$ and to $q_3$
enclose an angle close to $\pi / 2$.
Then Claim 2 and the arguments in part {\bf{F}}
of the proof of Lemma 2
yield a contradiction, for $D_d$ sufficiently large and $c_d > 0$ sufficiently
small.

Second, let $d=3$. We copy the proof of the case $d=2$.
Then either
$\ell = 2$ and $d - \ell + 2 = 3$,
or $\ell = 3$ and $d - \ell + 2 = 2$. Let the
$\sim $-equivalence classes on the primary colour classes of
$G_1^{*}{'}$
be represented either by the vertices
$\{ q_1,q_3,q_5 \} ,\,\,\{ q_2,q_4,q_6 \} $, or by the vertices
$\{ q_1,q_4 \}, \,\,\{ q_2,q_5 \}, \,\, \{ q_3,q_6 \} $, resp.
Here, all $q_j$-s belong to distinct primary colour
classes of $G_1^{*}{'}$,
which have $1,1,1,1,1,\lfloor {\text{const}} \cdot n \rfloor $
vertices, resp.
Up to notation, in both cases
we may assume that $q_6$ belongs to the last one of these primary colour
classes.
By Lemma 4, for $\ell = 2$, both
$\sphericalangle q_2 q_6 q_3 $ and $\sphericalangle q_4 q_6 q_3 $
are close to $ \pi /2 $, while for $\ell = 3$, both
$\sphericalangle q_3 q_6 q_4 $ and $\sphericalangle q_3 q_6 q_5 $
are close to $ \pi /2 $. Moreover, by Lemma 1 (3), $\ell = 2$
implies that
$\sphericalangle q_2 q_6 q_4 $ is close to $\pi /3$, while 
$\ell = 3$ implies that
$\sphericalangle q_4 q_6 q_5 $ is close to $\pi /3$.
Fix $q_1, \ldots , q_5$, and vary $q_6$ in its own primary
colour
class in $G_1^{*}{'}$, so that
it assumes $\lfloor {\text{const}} \cdot n \rfloor $
values (points). Then all these points lie
in the intersection of three spherical shells. These have
centres $q_2, q_4, q_3$ 
for $\ell = 2$, and centres $q_3, q_4, q_5$ for $\ell = 3$.
Moreover, their inner radii are some $t_{\kappa }$'s, and the
outer radii are the respective $(t_{\kappa } + c_dn^{1/d})$'s.
Further, for $\ell=2$ (and $3$), the angles enclosed by the
unit vectors pointing from $q_6$ to $q_2, q_3, q_4$ (and to
$q_3, q_4, q_5$, resp.,) are close to $\pi /3, \pi /2, \pi /2$.
Then Claim 2 and the arguments in part {\bf{F}}
of the proof of Lemma 2
yield a contradiction, for $D_d$ sufficiently large and
$c_d > 0$ sufficiently small.
\hfill
$\square $

\medskip


{\it{Proof of Theorem 1, continuation.}}
{\bf{Step 6}}.
By {\bf{Step 1}} (about tightness) and Lemmas 2, 7 and 8,
the proof of Theorem 1 for $d \ne 4,5$ follows.

\medskip


{\bf{Step 7}}.
Now we give the differences in the proof of Theorem 1
for the cases $d = 4,5$. Recall that the proof for Case I
already has been given in Lemma 2,
so we need to investigate Case II only (cf. {\bf{Step 3}}).

Analogously, as at the beginning of \S 3, we will have several
positive
constants, now depending on $d$ and $\varepsilon $, like
${\text{const}}_{d, \varepsilon }$, etc.
Of these, $C_{d, \varepsilon }$,
$D_{d, \varepsilon }$ will be fixed large constants, and 
$c_{d, \varepsilon }$ will be
a sufficiently small positive constant, in terms of the already
fixed values of all the other constants.

As in {\bf{Step 2}} of the proof of Theorem 1,
suppose, in order to obtain a contradiction, that the number of
edges of $G$ is greater than
$(n^2/2) (1 - 1/m(d - 1) + \varepsilon ) $. Then
by \cite
{1}, 
 Ch.\ 6, now $G$ contains a subgraph
$G_2 = K(\lfloor {\text{\rm const}}_{d, \varepsilon }
\cdot \log n \rfloor, \ldots, \lfloor
{\text{\rm const}}_{d, \varepsilon } \cdot \log n \rfloor )$, 
which is
a complete $ \left( m(d - 1) + 1 \right) $-partite graph, with
$\lfloor {\text{\rm const}}_{d, \varepsilon }
\cdot \log n \rfloor $
points in each colour class.
(For the dependence of ${\text{const}}_{d, \varepsilon }$
on $d$
and $\varepsilon $ in this statement,
the best known bound is given in \cite
{4}.) 

The primary (vertex) colouring and the secondary (edge)
colouring of $G_2$ are defined as for $G_1$ in {\bf{Step 2}}.
(Each edge of $G_2$ is coloured by $L$ or $R$.)
Analogously to the definition
of the subgraph $G_1^*$ of $G_1$ in {\bf{Step 2}}, for any
$j$, where $1 \le j \le m(d-1)+1$, we define an induced subgraph
$G_{2,j}^*$ of $G_2$ with the following properties.
Each primary
colour class of $G_2$, except the
$j$-th one, has exactly one point in $G_{2,j}^*$. Further, still
$\lfloor \text{\rm const}_{d, \varepsilon } \cdot \log n
\rfloor $
points of the $j$-th primary colour class of
$G_2$ belong to $G_{2,j}^*$. Moreover, the secondary colour
$L$ or
$R$ of an edge in $G_{2,j}^*$ depends only on the primary
colour classes the edge endpoints belong to.

Analogously to how we have defined the subgraph $G_1^{*}{'}$ of
$G_1^*$
in part {\bf{D}} of the proof of Lemma 2, now we define
the subgraph $(G^{*}_{2,j})'$ of $G_{2,j}^*$.
In what
follows, we will deal with the graphs $(G^{*}_{2,j}){'}$,
for $1 \le j \le m(d - 1) + 1$.

We want to show that 
for some $1 \le j \le m(d-1) + 1$ such a graph $(G_{2,j}^*)'$
cannot exist,
for $c_{d, \varepsilon } > 0$
a sufficiently small constant.
This contradiction will show that
the number of pairs $\{p_{i(1)}, p_{i(2)} \}$,
whose distances lie in the union of our
two intervals, is at most as large as stated in Theorem 1.

\medskip


{\bf{Lemma 9.}}
{\it{In Case II, for $d=4,5$, the upper estimate of Theorem 1 holds.}}

\medskip

{\it{Proof.}}
Suppose that none of the $\ell$ $\sim $-equivalence classes
of the primary colour classes of $(G_{2,j}^{*}){'}$
contains $d - \ell + 2$ primary colour classes (cf.\ Lemma 6).
Then we have $m(d - 1) + 1 \leq \ell(d - \ell +
1) \leq \lfloor \left( (d + 1)/2 \right)^2 \rfloor$. Thus,
by \thetag{1.2},
for $d
= 4$ we have $7 = m(3) + 1 \leq 6$, and for $d = 5$ we have
$11 = m(4) + 1
\leq 9$, i.e., in both cases we get a contradiction.
Therefore, {\it{one of the $\sim $-equivalence
classes of primary colour
classes, $C$, say, contains maximally many, i.e.,
$d - \ell + 2$ primary
colour
classes of $(G^{*}_{2,j}){'}$}}.

Let us choose from $C$ one primary colour class,
the $j(0)$-th one,
say, where $1 \le j(0) \le m(d-1) + 1$, and consider
$(G_{2,j(0)}^{*}){'}$.
We use the notation $r_1, \ldots , r_{\ell }$
introduced in {\bf{Step 5}}.

Let the $d - \ell + 2$ primary colour
classes in $C$
contain $d - \ell + 2$ points as follows. One of them is
$r_{j(0)}$,
and the others are
$q_{\ell + 1}, \ldots,  q_{d + 1}$.
Now consider the $d - \ell + 1$  vectors
$(q_j - r_{j(0)} )/d(q_j, r_{j(0)}) \in S^{d - 1}$,
for $j \in \{ \ell + 1, \ldots, d + 1 \} $,
and the $\ell - 1$ vectors $(r_{j(1)} - r_{j(0)})/
d(r_{j(1)}, r_{j(0)})
\in S^{d - 1}$, for $j(1) \in \{ 1, \ldots, \ell \}
\setminus \{ j(0) \} $. Let $M$ be the $d \times d$ matrix
formed by the above $(d - \ell + 1) + (\ell - 1) = d$
column vectors, in the above order. Then $|{\text{det}}\,M|
= [{\text{det}}\,(M'M)]^{1/2}$, where $M'$ is the transpose of
$M$. The entries of $M'M$ are the cosines of the angles formed
by the $d$ column vectors of $M$. The diagonal entries
of $M'M$ are equal to $1$. 
Outside the diagonal, in the intersection of the first
$d - \ell + 1$ rows and the first $d - \ell + 1$ columns,
as well as in the intersection of the last $\ell - 1$ rows
and the last $\ell - 1$ columns, by Lemma 1 (3),
the entries of $M'M$ are close to $1/2$. By Lemma 4,
the remaining entries of $M'M$ are close to $0$.

Let $N_0$ denote the $d \times d$ matrix, with the exact entries
$1, 1/2$ and $0$ at its respective positions.
Then ${\text{det}}\,(M'M)$ is close to
${\text{det}}\,N_0$, hence
$$
|{\text{det}}\,M| = [{\text{det}}\,(M'M)]^{1/2}
{\text { is close to }} [{\text{det}}\,(N_0)]^{1/2}\ (\in
[0, \infty ) ).
\tag 3.2
$$

Now we define a $d \times d$ matrix $M_0$ (it will not be unique)
as follows. Its first $d - \ell + 1$
column vectors are the edge vectors of a regular
$(d - \ell + 1)$-simplex of unit edge lengths in the coordinate subspace spanned by the
first $d - \ell + 1$ basic unit vectors, pointing from some  of
its vertices to all its other vertices. Moreover, its last
$\ell - 1$ column vectors are the edge vectors of a regular
$(\ell - 1)$-simplex of unit edge lengths in the
coordinate subspace spanned by the
last $\ell - 1$ basic unit vectors, pointing from some of
its vertices to all its other vertices. Then all these $d$
column
vectors form a base of ${\Bbb{R}}^d$, hence $|{\text{det}}\,
M_0|$ is some positive constant ${\text{const}}_{d, \ell }$,
independently of the choice of $M_0$. Since $\ell $ can assume
only finitely many values (cf.\ Lemma 5), therefore $|{\text{det}}\,M_0| \ge
{\text{const}}_d > 0$. Moreover, the entries of $M_0'M_0$ are
the cosines of the angles formed by the $d$ column vectors
of $M_0$. Hence, we have $M_0'M_0 = N_0$, which implies 
$$
[{\text{det}}\,(N_0)]^{1/2} = [{\text{det}}\,(M_0'M_0)]^{1/2} 
= |{\text{det}}\,M_0| \ge {\text{const}}\,_d > 0.
\tag 3.3
$$

By {\thetag{3.2}} and {\thetag{3.3}}, also
$$
|{\text{det}}\,M| \ge {\text{const}}_d > 0,
{\text{ provided }} D_{d, \varepsilon } {\text{ is}}
$$
$$
{\text{sufficiently large and }} c_{d, \varepsilon } > 0
{\text{ is sufficiently small}}.
$$

On the other hand, the
$\lfloor {\text{const}} _{d, \varepsilon }\cdot \log n
\rfloor $ points
of $(G^{*}_{2,j}){'}$ in its $j(0)$-th primary colour class
should be
contained in an intersection of $d$ spherical shells
(called $S_{j(h)}$ in part {\bf{E}} of the proof of Lemma 2).
These have centres
$q_{\ell + 1}, \ldots , q_{d + 1}$ and
$r_{j(1)}$ for
$j(1) \in \{ 1, \ldots , \ell \} \setminus \{ j(0) \}$,
inner radii some $t_{\kappa }$, and outer radii
(differently from part {\bf{E}} of the proof of Lemma 2)
the respective
$t_{\kappa } + c_{d, \varepsilon } (\log n)^{1/d}$.

Moreover,
the unit vectors pointing from $r_{j(0)}$ to
the above centres, are the column vectors of a $d \times d$
matrix, having a
determinant of absolute
value bounded from below by a positive number.
Then the slabs $S'_{j(h)}$ in Claim 2 will
be replaced by new slabs. More exactly, $r_{j(0)}$ replaces
$q_{{m(d - 1) + 1},1}$, the present $d$
centres
replace $q_{j(1)}, \ldots , q_{j(d)}$
in part {\bf{D}} of the proof of
Lemma 2, $\log n$ replaces $n$, and a suitable
half-space replaces $H^+$ (in part {\bf{E}} of the proof of
Lemma 2).
Then $\Pi $ in Claim 2 of the proof of Lemma 2 will
be replaced by a parallelepiped, circumscribed about a ball of
diameter $3c_{d, \varepsilon } (\log n)^{1/d}$.
Moreover, $\Pi '$ in part {\bf{F}} of the proof of Lemma 2 will
be replaced by a  parallelepiped, with
inradius
$3c_{d, \varepsilon }(\log n)^{1/d}/2 + 1/2$, hence
of volume at most 
$$
{\text{const}}_{d, \varepsilon } \cdot c_{d, \varepsilon } ^d
\log n .
$$

Thus, with these changes, the analogue of Claim 2 of the proof
of Lemma 2
(with the same proof as cited after Claim 2)
and the arguments in part {\bf{F}} of the proof
of Lemma 2
yield a contradiction. Namely, for $c_{d, \varepsilon } > 0$
sufficiently small, we have the following.
The parallelepiped replacing
$\Pi '$ has not enough volume in order to contain
$\lfloor {\text{const}}_{d, \varepsilon} \cdot \log n \rfloor $
disjoint open balls of unit diameter.
\hfill
$\square$

\medskip


{\it{Proof of Theorem 1, continuation.}}
{\bf{Step 8}}.
By {\bf{Step 1}} (about tightness) and Lemmas 2 and 9,
the proof of Theorem 1 for $d = 4,5$ follows.

Together with {\bf{Step 6}}, this completes
the proof of Theorem 1.
\hfill
$\square$

\medskip


\heading
\S4.  Proof of Theorem 2
\endheading

In this section, we present the proof of Theorem 2. The proof
falls into five simple steps marked as {\bf{Step 1}}, {\bf{Step 2}},
etc. 

{\it{Proof of Theorem 2.}}
{\bf{Step 1}}.
Recall that the tightness of Theorem 2 (A) and (B) was shown by
Constructions 2 and 3. It remains to establish that
$(d + 1)^k$ in (A) and $T\left( (d + 1)^k +
1, n \right) $ in (B) are upper bounds for the respective quantities.
For (B), this follows from (A),
by Tur\'an's theorem.

\medskip
{\bf{Step 2}}.
We need to show that, for $0 < \varepsilon < \varepsilon _{d,k}$,
where $\varepsilon _{d,k} > 0$ is sufficiently small, any
$(k, \varepsilon )$-distance set $P$ in ${\Bbb{R}}^d$
has a cardinality at most $(d + 1)^k$. We use induction on~$k$.

For $k = 1$, this statement is valid for $1 + \varepsilon < (1 + 2/
d) ^{1/2}$ (for $d$ even), or for $1 + \varepsilon <
[ 1 + 2(d + 2)/\left( d(d + 2) - 1 \right) ] ^{1/2}$ (for $d$ odd),
resp.\ (cf.\ Sch\"utte \cite
{24}, 
 Satz~3).

Now let $k \geq 2$.
We may suppose without loss of generality that $t_1 < \ldots < t_k$.
(If two of these numbers are equal, then the statement follows
by induction.)
We may and will suppose $\varepsilon \leq 1$.

Let $D_{d,k} > 0$ be a sufficiently large constant.
We distinguish two cases:
$$
{\text{\rm{Case I: }}} t_k/t_1 \leq D_{d,k},
$$
$$
{\text{\rm{Case II: }}} t_k/t_1 > D_{d,k}.
$$

\medskip

{\bf{Step 3}}.
In Case I, we prove


\medskip

{\bf{Lemma 10.}}
{\it{If $t_k/t_1 \leq D_{d,k}$, then the upper estimate stated
in {\bf{Step 2}} is valid.}}

\medskip


{\it{Proof.}}
The ratio of any two distances determined by $P$ is
at most $t_k(1+ \varepsilon )/t_1 \leq 2D_{d,k}$.
Then, for sufficiently small $\varepsilon > 0$, we get
$$
|P| \le m(d,k),
$$
by using the analogues of
the compactness considerations from the proof of Lemma 2,
Claim~1.
(Actually, only (i), the
analogue of (ii) with $1/(2D_{d,k})$, and (iii) from Claim 1
are needed.)

Further, by \thetag{1.1} we have
$$m(d,k) \leq {d + k\choose k} = \frac{d + 1}{1} \cdot
\ldots \cdot \frac{d + k}{k} = \left( \frac{d}{1} + 1 \right)
\cdot \ldots \cdot
\left( \frac{d}{k} + 1 \right) \leq (d + 1)^k.$$
Hence,
$$
|P| \leq m(d,k) \le (d + 1)^k ,
$$
as claimed in {\bf{Step 2}}.
\hfill
$\square$


\medskip

{\bf{Step 4}}.
In Case II, we prove


\medskip

{\bf{Lemma 11.}}
{\it{If $t_k/t_1 > D_{d,k}$, then the upper estimate from
{\bf{Step 2}} holds.}}


\medskip

{\it{Proof.}}
In Case II, there exists an integer $j \in \{ 1, \ldots ,k-1 \} $ such that
$t_{j + 1}/t_j > D_{d,k}^{1/(k - 1)}$.
We consider a colouring of the edges of the complete graph on
the vertex set $P$ with $k$ colours. Namely,
every edge $\{ p_{i(1)}, p_{i(2)} \} $
gets a colour $j$ with $d(p_{i(1)}, p_{i(2)}) \in [t_j, t_j(1
+ \varepsilon)]$. (Such a colouring is not necessarily unique,
but this makes no difference.)

Let us call a distance $d(p_{i(1)}, p_{i(2)})$ {\it{small}}
if its colour is at most
$j$, and {\it{large}} if its colour is at least $j + 1$.
The quotient of any large and any small distance is at least $t_{j +
1}/  \left( t_j(1 + \varepsilon) \right) > D_{d,k}^{1/(k - 1)} /(1
+ \varepsilon ) \geq D_{d,k}^{1/(k - 1)} /2 =: D_{d,k}'$, where
$D_{d,k}'$ is a large constant. In particular, we will assume
that $D'_{d,k}>1,$ which implies that $(0, t_j(1
+ \varepsilon)] \cap [t_{j + 1}, \infty ) = \emptyset $. Thus,
the length
$d(p_{i(1)}, p_{i(2)})$ uniquely determines whether it is a small or a large
distance.
From now on, we also assume that $D'_{d,k}>2.$ This yields that
{\it{every large distance is more than twice as large as
every small distance}}.

This implies that we can define an equivalence
relation $\sim$ on the points $p_{i(1)}$, $p_{i(2)} \in P$.

\medskip


{\bf{Definition 5.}}
For $p_{i(1)}, p_{i(2)} \in P$ we write
$p_{i(1)} \sim p_{i(2)}$ if either $i(1) = i(2)$, or $d(p_{i(1)},
p_{i(2)})$
is a small distance. By the last italicized text,
$\sim $ is an equivalence relation on $P$.

\medskip


In each $\sim$-equivalence class of the points $p_i \in P$, each edge has a colour
at most $j$. Thus, each $\sim$-equivalence class is a
$(j, \varepsilon )$-distance set.
Since $j \le k-1$, by the induction hypothesis we have, for $\varepsilon > 0$
sufficiently small, that the cardinality of any $\sim$-equivalence class is at most
$(d + 1)^j$.

Now consider a set of representatives from each 
$\sim$-equivalence class.
In this set, each edge has a colour
at least $j + 1$, so it is a
$(k - j, \varepsilon )$-distance set.
Since $k - j \le k-1$, by the induction hypothesis we have, for $\varepsilon >
0$ sufficiently small, that the cardinality of this set is at
most $(d + 1)^{k - j}$.

Using the results of the last two paragraphs, we obtain the 
following. For $\varepsilon > 0$
sufficiently small, $|P|$
is at most the number of $\sim$-equivalence
classes times the maximum cardinality of a $\sim$-equivalence class. That is,
$$
|P| \leq (d + 1)^{k - j}(d + 1)^j = (d + 1)^k,
$$
as asserted in {\bf{Step 2}}.
\hfill
$\square$
\medskip


{\bf{Step 5}}. Now Theorem 2 follows from {\bf{Steps 1}},
{\bf{2}}, and Lemmas 10 and 11.
\hfill $\square$


\medskip
\heading
\S5. Concluding remarks
\endheading

{\bf{1.}} Suppose that neither $k$
nor $d$ is much larger than the
other.
It seems likely that in this case 
one can obtain reasonably good constructions for $Q$
in Construction 1 the following way. 
Suppose that $d = d(1) + \ldots + d(h)$ and $k = k(1) + \ldots + k(h)$, where
all $d(g)$ and $k(g)$, for $1 \le g \le h$, are natural numbers. Then ${\Bbb{R}}^d =
{\Bbb{R}}^{d(1) + \ldots + d(h)} = {\Bbb{R}}^{d(1)} \oplus \ldots
\oplus {\Bbb{R}}^{d(h)}$.
In each ${\Bbb{R}}^{d(g)}$, for $1 \le g \le h$,
we take a subset $Q_g$. Here, each $Q_g$
is one of the examples from Construction 1$'$
or Construction 1$''$, with all distances in $Q_g$ lying 
in the union of $k(g)$ intervals of arbitrarily small lengths. We scale
$Q_1, \ldots , Q_h$ in such a way that for each $1 \le g \le h-1$, the
maximal distance in $Q_g$ is much smaller than the minimal distance in
$Q_{g+1}$.
Moreover, all distances in $Q_g$
still belong to the union of $k(g)$ intervals of arbitrarily small lengths.
Let $Q := \oplus _{g=1}^h Q_g$. For any two
distinct points $q(1) = \oplus _{g=1}^h q_{g(1)}, \,\,
q(2) = \oplus _{g=1}^h q_{g(2)} \in Q$, there is a largest
$g \in \{ 1, \ldots , h \} $ such that $q_{g(1)} \ne q_{g(2)}$.
Then the
distance between $q(1)$ and $q(2)$ is arbitrarily close to the distance between
$q_{g(1)}$ and $q_{g(2)}$. Therefore, all distances determined by $Q$ lie in the
union of $k(1) + \ldots +k(h) = k$ intervals of arbitrarily small length.

\medskip

{\bf{2.}}
A related question was studied by Pach, Radoi\v ci\'c and
Vondr\'ak \cite
{22}, 
 \cite
{23}. 
 They proved that for any $d \geq 2$ and any
$0 < \gamma < 1/4$, the
following statement holds. Suppose that
in an $n$-element separated point set $P \subset
{\Bbb{R}}^d$ there are at least $ \gamma n^2$
point pairs whose distances differ by at most $1$. Then the
diameter of $P$ is at least $\text{\rm const}_{d, \gamma }
\cdot n^{2/(d - 1)}$. Apart from the value of the constant, this
bound is tight for all $d \geq 2$ and all $0 < \gamma < 1/4$.

\medskip


{\bf{3.}} Another related question is treated in \cite
{11}. 
 Suppose that in a separated $n$-element
point set $P$ in the plane, the number of pairs that
determine a distance
nearly equal to one of $t_1 < \ldots < t_k$ is maximal.
Does it follow that
then we have, ``approximately,'' $t_2 = 2t_1, \ldots, t_k
= kt_1$ (as in the example after Theorem B)?
In this direction, they proved the following.
Let $\delta > 0$ and suppose that for any $1 \le i(1) \leq i(2) <
i(3) \le k$, the inequality $| t_{i(3)} / (t_{i(1)} +
t_{i(2)}) - 1 | > \delta$ holds.
Then, for $n \ge n_{k, \delta } $,
the number of unordered pairs that determine a distance belonging to $[t_1, t_1 +
1] \cup \ldots \cup [t_k, t_k + 1]$, is at most $n^2/4 + \text{\rm
const}_{k,\delta} \cdot n$. This bound is sharp, up to the value of
$\text{\rm const}_{k, \delta} > 0$.
It is easy to see that if $t_{i(3)} = t_{i(1)} + t_{i(2)}$ holds for some
$i(1) \leq i(2) < i(3)$, then the
number of pairs with the above property can attain
$\lfloor n^2/3\rfloor$.
\medskip


{\bf{4.}} We pose the following

{\bf{Question}}. What would be the
results analogous
to our Theorem 2, for unions of intervals of the form
$[t_1, t_1^{1 + \varepsilon }] \cup \ldots \cup
[t_k, t_k^{1 + \varepsilon }]$, for $\varepsilon > 0$?


{\bf{Acknowledgements.}} This research has been partially supported
over the years by several OTKA (Hungarian National Foundation
for Scientific
Research) and NKFIH (National Office of Research, Development,
and
Innovation) grants, most recently, by grant no. K131529, and by ERC
Advanced Grant "GeoScape" no. 882971.


\head
References
\endhead

\frenchspacing

\parindent=45truept

\item{\hbox to41truept{[1]
\hfill}} {B. Bollob\'as,}
{\it Extremal graph theory,} London Math. Soc. Monographs {\bf{11}},
Academic Press, Inc.
[Hartcourt Brace Jovanovich, Publishers], London-New York, 1978,
MR~{\bf{80a:}}05120

\item{\hbox to41truept{[2] 
\hfill}} {A. Blokhuis,}
A new upper bound for the cardinality of 2-distance sets in Euclidean space,
in:
{\it Convexity and Graph Theory (Jerusalem, 1981)},
North-Holland Math. Stud. {\bf 87}, {\it Ann. Discr. Math.} {\bf 20},
North-Holland, Amsterdam, 1984, 65--66, MR~{\bf{86h:}}52008

\item{\hbox to41truept{[3] 
\hfill}} {A. Blokhuis, J. J. Seidel,}
Few-distance sets in ${\bold{R}}^{p,q}$,
in: {\it Symposia Mathematica}, {\bf{38}} (Rome, 1983),
Academic Press, London - New York, 1986, 145--158, MR~{\bf{88d:}}52011

\item{\hbox to41truept{[4] 
\hfill}} {V. Chv\'atal, E. Szemer\'edi,}
On the Erd\H os-Stone theorem, {\it J. London Math. Soc.} {\bf 23}
(1981), 207-214, MR~{\bf{82f:}}05057

\item{\hbox to41truept{[5] 
\hfill}} {H. T. Croft}, 9-point and 7-point
configurations in 3-space,
{\it Proc. London Math. Soc. (3)} {\bf 12} (1962), 400--424, corr.: ibid. {\bf
13} (1963), 384, MR~{\bf{27{\#}}}
\newline
5167

\item{\hbox to41truept{[6] 
\hfill}} {P. Erd\H{o}s,} On a lemma of Littlewood and Offord,
{\it Bull. Amer. Math. Soc.} {\bf 51} (1945), 898--902, MR~{{\bf 7,}309j}

\item{\hbox to41truept{[7] 
\hfill}} {P. Erd\H{o}s,} On sets of distances of $n$
points, {\it Amer. Math. Monthly}\,\, {\bf 53} (1946), 248-250, MR~{\bf{7,}}741c

\item{\hbox to41truept{[8] 
\hfill}} {P. Erd\H{o}s, L. M. Kelly,}
Elementary Problems and Solutions: Solutions: E735. Isosceles $n$-points,
{\it Amer. Math. Monthly}\,\, {\bf 54} (1947), 227--229, MR~{\bf{1526679}}

\item{\hbox to41truept{[9] 
\hfill}} {P. Erd\H{o}s, E. Makai, Jr., J. Pach,}
Nearly equal distances in the plane,
{\it Combinatorics, Probability and Computing} {\bf 2} (1993), 401--408,
MR~{\bf{95i:}}52018

\item{\hbox to41truept{[10] 
\hfill}} {P. Erd\H{o}s, E. Makai, Jr., J. Pach,}
Two nearly equal distances in ${\Bbb{R}}^d$, arXiv:
\newline
1901.01055

\item{\hbox to41truept{[11] 
\hfill}} {P. Erd\H{o}s,
E. Makai, Jr., J. Pach,}
Nearly equal distances in the plane II, arXiv:2112.08852

\item{\hbox to41truept{[12] 
\hfill}} {P. Erd\H{o}s, E. Makai, Jr., J. Pach,
J. Spencer,}
Gaps in difference sets and the graph of nearly equal distances,
{\it Applied Geometry and Discr. Math., The V. Klee Festschrift, DIMACS Series
in Discr. Math. and Theoretical Computer Science}
 {\bf 4} (1991), 265--273, MR~{\bf{92i:}}52021

\item{\hbox to41truept{[13] 
\hfill}}P. Erd\H os, A. H. Stone, On the structure of linear graphs, {\it Bull. Amer. Math. Soc.} {\bf 52} (1946), 1087--1091, MR~{{\bf 8,}333b}

\item{\hbox to41truept{[14] 
\hfill}} {N. Frankl, A. Kupavskii,} Nearly $k$-distance sets, 
Acta Math. Univ. Comenian. (N. S.) {\bf{88}} (2019) (3), 689-693,
MR~{\bf{4012869}}

\item{\hbox to41truept{[15] 
\hfill}} {N. Frankl, A. Kupavskii,} Nearly $k$-distance sets, 
Discrete Comput. Geom. {\bf{70}} (2023) (3), 455-494, 
MR~{\bf{4650014}}

\item{\hbox to41truept{[16] 
\hfill}} {L. Guth, N. H. Katz,}
On the Erd\H os
distinct distances problem in the plane, Ann. of Math. (2) {\bf{181}} (2015),
155-190, MR~{\bf{3272924}}

\item{\hbox to41truept{[17] 
\hfill}} {J. W. Hurewicz, H. Wallman},
{\it Dimension Theory}, Princeton Math. Series {\bf{4}}, Princeton University Press,
Princeton, N. J., 1941, MR~{\bf{3,}}312b

\item{\hbox to41truept{[18] 
\hfill}} {P. Lison\v ek,}
New maximal two-distance sets,
{\it J. Combin. Th., Ser. A} {\bf 77} (1997), 318--338, MR~{\bf{98a:}}51014

\item{\hbox to41truept{[19] 
\hfill}} J. E. Littlewood, A. Offord, On the number of real roots of a random algebraic equation. III. {\it Rec. Math. [Mat. Sbornik] N.S.} {\bf 12} (1943), 277--286, MR~{\bf{5,}}179h

\item{\hbox to41truept{[20] 
\hfill}} {E. Lutwak,}
Selected affine isoperimetric inequalities, {\it Handbook of Convex
Geometry} (Eds. P.M. Gruber, J. M. Wills), North Holland,
Amsterdam etc., 1993, Vol. A, Chapter 1.5, 151-176,
MR~{\bf{94h:}}52014

\item{\hbox to41truept{[21] 
\hfill}} {E. Makai, Jr., J. Pach, J. Spencer,}
New results on the distribution of distances determined by separated point
sets,
{\it Paul Erd\H{o}s and His Math. II}, Bolyai Soc. Math. Studies {\bf{11}}.
(Ed. by G. Hal\'asz, L. Lov\'asz, M. Simonovits, V. T. S\'os.)
Springer, Berlin etc., J. Bolyai Math. Soc., Budapest, 2002, 499--511,
MR~{\bf{2004a:}}52026

\item{\hbox to41truept{[22] 
\hfill}} {J. Pach, R. Radoi\v ci\'c, J. Vondr\'ak,}
Nearly equal distances and Szemer\'edi's regularity lemma,
{\it Comput. Geom.} {\bf 34} (2006), 11--19, MR~{\bf{2006i:}}52022

\item{\hbox to41truept{[23] 
\hfill}} {J. Pach, R. Radoi\v ci\'c, J. Vondr\'ak,}
On the diameter of separated point sets with many equal distances,
{\it European J. Combin.} {\bf 27} (2006), 1321--1332. MR~{\bf{2007f:}}51026

\item{\hbox to41truept{[24] 
\hfill}} {K. Sch\"utte,}
Minimale Durchmesser endlicher Punktmengen mit vor\-ge-
\newline
schrie\-be\-nem Mindestabstand 
(The minimal diameter of finite point sets with prescribed
minimal distance, in German)
{\it Math. Ann.} {\bf 150} (1963), 91--98, MR~{\bf{26{\#}}}5479

\item{\hbox to41truept{[25] 
\hfill}} {J. Solymosi, V. H. Vu,}
Near optimal bounds for the Erd\H os distinct distances problem in high
dimensions,
{\it Combinatorica} {\bf 28} (2008), 113--125, MR~{\bf{2009f:}}52042

\item{\hbox to41truept{[26] 
\hfill}}  P. Tur\'an, On an extremum problem
in graph theory
(in Hungarian),
{\it Ma\-te\-ma\-ti\-kai \'es Fizikai Lapok} {\bf 48} (1941), 436--452,
MR~{\bf{8,}}284j

\end

\enddocument